\documentclass[11pt,reqno,a4paper]{amsart}

\usepackage{amsthm,amsmath,amsfonts,amssymb,epsfig,graphicx,latexsym,float,color}
\usepackage{mathtools,multicol} 

\usepackage{enumerate,relsize,mathrsfs}
\usepackage{sidecap}
\usepackage{multirow}

\newtheorem{definition}{Definition}[section]

\newtheorem{remark}[definition]{Remark}
\newtheorem{example}[definition]{Example}
\newtheorem{notation}[definition]{Notation}

\usepackage{bbold}

\newcommand{\mdef}{\ensuremath{\mathbin{:=}}}

\newcommand{\dx}[1]{\ensuremath{\,\mathrm{d}{#1}}}
\newcommand{\pp}{\ensuremath{\mathbb{P}}}
\newcommand{\ii}{\ensuremath{\mathbb{I}}}
\newcommand{\nn}{\ensuremath{\mathbb{N}}}
\newcommand{\rr}{\ensuremath{\mathbb{R}}}

\renewcommand{\vec}[1]{\ensuremath{\mathbf{#1}}}
\newcommand{\tens}[1]{\ensuremath{\mathcal{#1}}}
\newcommand{\mat}[1]{\ensuremath{\mathbb{#1}}}

\newcommand{\transpose}{\ensuremath{\mathsf{T}}}

\newcommand{\STAB}[1]{\begin{tabular}{@{}c@{}}#1\end{tabular}}

\title{Extended Group Finite Element Method}
\author[K.~Tolle]{Kevin Tolle$^{\dag,1}$}
\author[N.~Marheineke]{Nicole Marheineke$^1$}
\date{\today\\%
$^\dag$ \textit{Corresponding author}, email: tolle@uni-trier.de, phone: +49\,651\,201\,3472\\%
$^1$ Universit\"at Trier, FB IV - Mathematik, Lehrstuhl Modellierung und Numerik, D-54286 Trier, Germany%
}

\begin{document}

\begin{abstract}
    Interpolation methods for nonlinear finite element discretizations are commonly used to eliminate the computational costs associated with the repeated assembly of the nonlinear systems. While the group finite element formulation interpolates nonlinear terms onto the finite element approximation space, we propose the use of a separate approximation space that is tailored to the nonlinearity. In many cases, this allows for the exact reformulation of the discrete nonlinear problem into a quadratic problem with algebraic constraints. Furthermore, the substitution of the nonlinear terms often shifts general nonlinear forms into trilinear forms, which can easily be described by third-order tensors. The numerical benefits as well as the advantages in comparison to the original group finite element method are studied using a wide variety of academic benchmark problems. 
\end{abstract}

\maketitle

\noindent
\textsc{AMS-Classification:} 65N30, 65Y20  
 

\section{Introduction}

A wide variety of mathematical models use partial differential equations (PDEs) to describe physical processes. A powerful tool for the numerical solution of PDEs is the finite element method, see, for example, \cite{Ern2004,Brenner2008} for theoretical aspects and \cite{Logg2012} for practical implementation aspects. While linear problems are straightforward and can be solved directly, nonlinear problems present more difficulty. For example, unlike their linear counterparts, nonlinear problems require iterative methods, such as fixpoint iterations or Newton's method. For finite element discretizations, this means that the computationally expensive assembly of the system is required for each iteration. There are many different approaches to improve the efficiency of the assembly process such as parallelization and optimized algorithms using form compilers, see \cite{Logg2012}. However, interpolation methods offer a different approach by circumventing the assembly of the nonlinear system in each iteration.

\begin{figure}[!tb] \centering
	\includegraphics{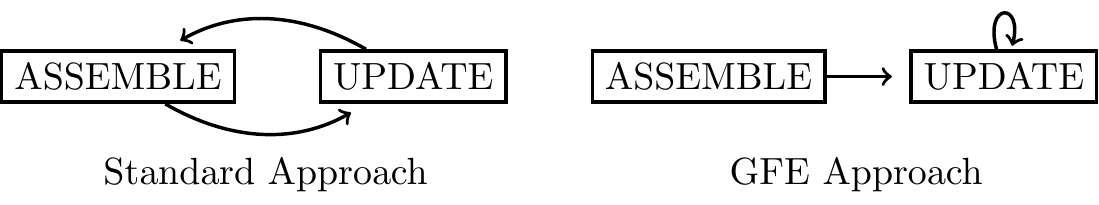}

	\caption{Simple sketch of the workflow for the standard Galerkin approach (SGA) using finite elements and the group finite element method (GFEM).}
\label{fig:flowchart}
\end{figure}

The group finite element formulation, also known as the product approximation technique, is an interpolation method that eliminates the repeated assembly of the nonlinear system, see \cite{Christie1981,Fletcher1983} and references therein. Although there is no direct connection, the group finite element method (GFEM) can be seen as a specific application of the interpolation procedure described in \cite{Douglas1975}. The GFEM replaces groups of variables with an interpolant on the finite element approximation space, which allows the system to be precomputed. While the standard Galerkin approach (SGA) reassembles the system after each update to the approximate solution, the GFEM assembles the interpolated system once before the iterative method is started, see Figure \ref{fig:flowchart}. This approach is still an active and relevant research topic. For example, it has been combined with proper orthogonal decomposition in the context of nonlinear model order reduction in \cite{Dickinson2010}. Additionally, a convergence analysis of the group finite element formulation is presented in \cite{Barrenechea2017}, which also briefly mentions the use of the GFEM in implicit flux-corrected transport schemes, see, for example, \cite{Kuzmin2012}.

This paper investigates the benefits of extending the group finite element formulation to allow for general finite element approximation spaces for the nonlinear terms. To the best of our knowledge, previous work has only focused on other approaches: For example, a variety of interpolation techniques, such as least-squares fitting, Lagrange or Hermitian interpolation and cubic spline interpolation, for nonlinear functions are investigated in \cite{Douglas1975}. Alternatively, the GFEM interpolates nonlinear terms onto the underlying finite element approximation space, see \cite{Christie1981,Fletcher1983}. Our generalization, which we call the extended group finite element method (EGFEM), allows for more flexibility and broadens the scope of problems that can be treated. Beyond the decrease in the computational effort by eliminating repeated assembly, no sources known to us discuss the unique structure of the interpolated formulation with respect to multi-linear forms. Therefore, an additional aspect that we study is the underlying tensor structure that arises by introducing new variables in place of the nonlinearities. This structure plays an important role in the efficiency of our method, while also showing great promise for model order reduction techniques. 

This paper is organized as follows. The mathematical setting and the original group finite element formulation are described in Section \ref{sec:background}. We present our extended group finite element method in Section \ref{sec:egfem} with a focus on the construction and structure of the general approximations as well as the vectorized evaluation of the nonlinear terms. The numerical realizations of our extended group finite element method based on Picard iterations and Newton's method are compared in Section \ref{sec:numeric}. Finally, Section \ref{sec:results} contains a selection of numerical experiments, ranging from the Burgers' equation, which serves as a benchmark problem for the group finite element method, to the $ p $-Laplace equation, which belongs to a class of problems that the original GFEM struggles with. At the end of Section \ref{sec:results}, we present some guidelines based on the numerical experiments.

\section{Background} \label{sec:background}

We present a model problem, which represents a class of nonlinear PDEs, that we consider in this work in Section \ref{subsec:problem}. For simplicity, our considerations are restricted to stationary boundary-value problems, time-dependent and/or more general nonlinear problems can be handled similarly. The model problem is discretized using the finite element method, which leads to a nonlinear system that is solved iteratively. Because each iteration requires the costly assembly of the system, we briefly present the group finite element method in Section \ref{subsec:gfem}. This method, which is also known as the product approximation, interpolates nonlinear terms onto the finite element approximation space so that assembly only needs to be performed once.

\subsection{Problem Description} \label{subsec:problem}

Let $ \Omega \subset \rr^{n} $ with $ n \in \nn $ be a bounded domain with sufficiently smooth boundary $ \Gamma = \Gamma_{\mathrm{D}} \cup \Gamma_{\mathrm{N}} $ such that $ \Gamma_{\mathrm{D}} \cap \Gamma_{\mathrm{N}} = \emptyset $. We consider a nonlinear partial differential equation of the form
\begin{subequations} \label{eqn:nonlinear_pde}
\begin{align}
    -\nabla_{x} \cdot \left(a(x,u(x),\nabla u(x)) \, \nabla_{x} b(u(x))\right) + c(x,u(x),\nabla u(x)) &= d(x), & x &\in \Omega,
\intertext{with mixed boundary conditions}
    u(x) &= u_{\mathrm{D}}(x), & x &\in \Gamma_{\mathrm{D}}, \\
    a(x,u(x),\nabla u(x)) \nabla_{x} b(u(x)) \cdot \vec{n} + g(x,u(x)) &= h(x), & x &\in \Gamma_{\mathrm{N}},
\end{align}
\end{subequations}
where $ \vec{n} $ denotes the unit normal vector along the boundary. For better readability, we omit the argument $ x \in \Omega $ in the following. We assume that the scalar-valued functions in \eqref{eqn:nonlinear_pde} ensure the existence of a weak solution $ u $ in the trial space $ V $. A weak solution $ u \in V $ fulfills $ u|_{\Gamma_{\mathrm{D}}} = u_{\mathrm{D}} $ as well as the variational problem:
\begin{gather} \label{eqn:weak_form}
    \int_{\Omega} a(u,\nabla u) \, \nabla_{x} b(u) \cdot \nabla v + \left(c(u,\nabla u) - d\right) v \dx{x} = \int_{\Gamma_{\mathrm{N}}} \hspace{-0.5em} \left(h - g(u)\right) v \dx{x}
\end{gather}
for all $ v \in V_{0} \mdef \{v \in V : v|_{\Gamma_{\mathrm{D}}} = 0\} $. While the Dirichlet boundary conditions are enforced separately, the Neumann/Robin boundary conditions are naturally embedded in the weak formulation \eqref{eqn:weak_form}. For more details on the theoretical aspects of nonlinear partial differential equations, we refer the reader to \cite{Roubicek2013,Evans2010} and references therein.

The spatial discretization of \eqref{eqn:weak_form} is achieved by restricting the infinite-dimensional space $ V $ to a discrete (finite-dimensional) subspace $ V_{\mathrm{h}} \subset V $. The finite element method constructs this discrete space by introducing a partitioning $ \Omega_{\mathrm{h}} $ of the domain and using local function spaces on the $n$-dimensional simplices, see, e.g., \mbox{\cite{Ern2004,Logg2012}} for details. For the resulting discrete space, we assume that a set of nodal basis functions $ \{\phi_{i}\}_{i=1}^{N_{\mathrm{u}}} $ is given, such that $ V_{\mathrm{h}} = \mathop{\mathrm{span}} \{\phi_{1},\ldots,\phi_{N_{\mathrm{u}}}\} $. The approximate solution $ u_{\mathrm{h}} \in V_{\mathrm{h}} $ is represented as a linear combination of the basis functions, i.e.,
\begin{gather*}
    u(x) \approx u_{\mathrm{h}}(x) = \sum_{i=1}^{N_{\mathrm{u}}} u_{i} \, \phi_{i}(x).
\end{gather*}
Without loss of generality let the first $ N_{\mathrm{D}} < N_{\mathrm{u}} $ basis functions be ordered so that they belong to nodes $ x_{i} \in \Omega_{\mathrm{h}} $ on the discrete Dirichlet boundary $ \Gamma_{\mathrm{D},\mathrm{h}} $, which consists of the interfaces, e.g., edges in two-dimensions and faces in three-dimensions, that approximate the Dirichlet boundary. Thus, we can introduce the discrete test space $ V_{\mathrm{h},0} \mdef \mathop{\mathrm{span}} \{\phi_{N_{\mathrm{D}}+1},\ldots,\phi_{N_{\mathrm{u}}}\} $. The approximate solution $ u_{\mathrm{h}} \in V_{\mathrm{h}} $ fulfills $ u_{i} = u_{\mathrm{D}}(x_{i}) $ for $ 1 \leq i \leq N_{\mathrm{D}} $ as well as
\begin{gather} \label{eqn:discrete_problem}
    \int_{\Omega_{\mathrm{h}}} a(u_{\mathrm{h}},\nabla u_{\mathrm{h}}) \, \nabla_{x} b(u_{\mathrm{h}}) \cdot \nabla \phi_{i} + \left(c(u_{\mathrm{h}},\nabla u_{\mathrm{h}}) - d\right) \phi_{i} \dx{x} = \int_{\Gamma_{\mathrm{N},\mathrm{h}}} \hspace{-1.0em}  \left(h - g(u_{\mathrm{h}})\right) \phi_{i} \dx{x}
\end{gather}
for all $ \phi_{i} \in V_{\mathrm{h},0} $, where the discrete Neumann/Robin boundary $ \Gamma_{\mathrm{N},\mathrm{h}} $ is defined similarly to $ \Gamma_{\mathrm{D},\mathrm{h}} $. The integrals in \eqref{eqn:discrete_problem} are straightforwardly defined as the sum of all the integrals over the individual elements and interfaces, respectively. These nonlinear equations are solved numerically using some form of iterative method, which may either be based on a Picard/fix-point or Newton iteration, see \cite{Logg2012} for more details. However, the major difficulty for solving these nonlinear finite element problems lies in the repeated computation of the nonlinear terms, i.e.,
\begin{gather} \label{eqn:nonlinear_terms}
\begin{aligned}
    \left(\vec{K}(a,b,\vec{u})\right)_{i} &\mdef \sum_{K \in \Omega_{\mathrm{h}}} \int_{K} a(\sum_{j = 1}^{N_{\mathrm{u}}} u_{j} \, \phi_{j},\sum_{j = 1}^{N_{\mathrm{u}}} u_{j} \, \nabla \phi_{j}) \, \nabla_{x} b(\sum_{j = 1}^{N_{\mathrm{u}}} u_{j} \, \phi_{j}) \cdot \nabla \phi_{i} \dx{x}, \\ 
    \left(\vec{M}(c,\vec{u})\right)_{i} &\mdef \sum_{K \in \Omega_{\mathrm{h}}} \int_{K} c(\sum_{j = 1}^{N_{\mathrm{u}}} u_{j} \, \phi_{j},\sum_{j = 1}^{N_{\mathrm{u}}} u_{j} \, \nabla \phi_{j}) \, \phi_{i} \dx{x}, \\
    \left(\vec{g}(\vec{u})\right)_{i} &\mdef \sum_{S \in \Gamma_{\mathrm{N},\mathrm{h}}} \int_{S} g(\sum_{j = 1}^{N_{\mathrm{u}}} u_{j} \, \phi_{j}) \, \phi_{i} \dx{x},
\end{aligned}
\end{gather}
where $ \vec{u} $ denotes $ (u_{1},\ldots,u_{N_{\mathrm{u}}})^{\transpose} $. Because these terms often require numerical integration in order to evaluate them, alternative methods that eliminate the repeated assembly of the discrete system are of great importance.

\subsection{Group Finite Element Method} \label{subsec:gfem}

The group finite element method simplifies the nonlinear problems in order to avoid the assembly of the nonlinear system in each iteration of the numerical solution method, see \cite{Christie1981,Fletcher1983}. Assume that the given basis $ \{\phi_{i}\}_{i=1}^{N_{\mathrm{u}}} $ is a nodal basis, i.e., $ \phi_{i}(x_{j}) = \delta_{ij} $ at the degrees of freedom $ x_{j} $ associated to the finite element ansatz. This allows \emph{groups} (products) of variables to be easily interpolated onto the trial space $ V_{\mathrm{h}} $, e.g., the scalar-valued nonlinear function $ f(u_{\mathrm{h}}) $ is approximated through
\begin{gather*}
    f(u_{\mathrm{h}}(x)) = f(\sum_{i=1}^{N_{\mathrm{u}}} u_{i} \, \phi_{i}(x)) \approx \sum_{i=1}^{N_{\mathrm{u}}} f(u_{i}) \, \phi_{i}(x),
\end{gather*}
which follows from the fact that $ f(u_{\mathrm{h}}(x_{i})) = f(u_{i}) $. The advantage of this method in comparison to the standard Galerkin approach (SGA) is that the nonlinear terms become (multi-)linear forms that can be precomputed, e.g.,
\begin{gather*}
    \int_{\Omega_{\mathrm{h}}} f(u_{\mathrm{h}}) \, \phi_{i} \dx{x} \approx \int_{\Omega_{\mathrm{h}}} \left(\sum_{j=1}^{N_{\mathrm{u}}} f(u_{j}) \, \phi_{j}\right) \phi_{i} \dx{x} = \sum_{j=1}^{N_{\mathrm{u}}} f(u_{j}) \left(\int_{\Omega_{\mathrm{h}}} \phi_{j} \, \phi_{i} \dx{x} \right).
\end{gather*}
Because the integration only needs to be performed once, the computational overhead of each iteration is reduced noticeably. Additionally, the nonlinearity is embedded into the point-wise evaluation of $ f(u_{\mathrm{h}}) $ at the nodes.

Although the GFEM is approximative in nature, there are select examples that display greater accuracy in the nodal values, when compared to the standard approach, see \cite{Christie1981,Fletcher1983}. However, these results depend directly on the norm used to compute the error. For one of the examples, the $ L^{2} $-norm converges up to two orders of magnitude slower than the standard Galerkin approach, as shown in \cite{Christie1981}. Our goal is to extend the GFEM in such a way that the reformulated problem is as exact as possible, while maintaining the reduced complexity and performance gains of the original method. Additionally, we generalize the form of the problems considered in \cite{Christie1981} to include the more general form \eqref{eqn:nonlinear_pde}, which leads to trilinear forms.

\section{Extended Group Finite Element Method} \label{sec:egfem}

In this section, we extend the group finite element method in order to allow for more generality in the nonlinear terms and their approximations. The latter is achieved by introducing an approximation space $ W_{\mathrm{h}} $ for the nonlinearity that is tailored specifically to the problem instead of recycling the trial space $ V_{\mathrm{h}} $. We discuss the implications and handle the selection of this general approximation in detail for a scalar-valued nonlinear function $ f $ that acts on the solution $ u $ and/or its gradient $ \nabla u $ in Section \ref{subsec:egfem_choice}. Extending our approach to functions with an explicit dependence on the position $ x \in \Omega_{h} $ is straightforward. Afterwards, we briefly present an efficient approach for the point-wise evaluation of $ f(u_{\mathrm{h}},\nabla u_{\mathrm{h}}) $ with the help of interpolation operators in Section \ref{subsec:algebraic}. By considering a general model problem, multi-linear forms result after applying the extended group finite element method. In our case, we attain trilinear forms, which are represented through third-order tensors. The resulting tensor formulation is developed and compared to the standard Galerkin approach in Section \ref{subsec:egfem_tensors}. 

\subsection{Choice of Ansatz Spaces for Nonlinear Terms} \label{subsec:egfem_choice}

The group finite element method interpolates nonlinear functions onto the trial space $ V_{\mathrm{h}} $. While this simple choice may be viable in many cases, we introduce a separate approximation space $ W_{\mathrm{h}} $ that is spanned by the nodal basis $ \{\eta_{k}^{\mathrm{f}}\}_{k=1}^{N_{\mathrm{f}}} $ with the degrees of freedom $ \{x_{k}^{\mathrm{f}}\}_{k=1}^{N_{\mathrm{f}}} $. Interpolating the nonlinear function $ f(u_{\mathrm{h}},\nabla u_{\mathrm{h}}) $ onto $ W_{\mathrm{h}} $ instead of $ V_{\mathrm{h}} $ delivers the following approximation:
\begin{gather*}
    f(u_{\mathrm{h}}(x),\nabla u_{\mathrm{h}}(x)) = f(\sum_{i=1}^{N_{\mathrm{u}}} u_{i} \, \phi_{i}(x),\sum_{i=1}^{N_{\mathrm{u}}} u_{i} \, \nabla \phi_{i}(x)) \approx \sum_{k=1}^{N_{\mathrm{f}}} f_{k} \, \eta_{k}^{\mathrm{f}}(x) = f_{\mathrm{h}}(x).
\end{gather*}
The coefficients $ \vec{f} = (f_{1},\ldots,f_{N_{\mathrm{f}}})^{\transpose} $ are defined through the point-wise evaluation of the function $ f $ on the degrees of freedom, i.e.,
\begin{gather} \label{eqn:egfem_coeffs}
    f_{k} \mdef f(u_{\mathrm{h}}(x_{k}^{\mathrm{f}}),\nabla u_{\mathrm{h}}(x_{k}^{\mathrm{f}})) = f(\sum_{i=1}^{N_{\mathrm{u}}} u_{i} \, \phi_{i}(x_{k}^{\mathrm{f}}),\sum_{i=1}^{N_{\mathrm{u}}} u_{i} \, \nabla \phi_{i}(x_{k}^{\mathrm{f}})),
\end{gather}
which is motivated by the fact that $ f_{\mathrm{h}}(x_{k}^{\mathrm{f}}) = f_{k} $ by construction.  A nice feature of our extended group finite element method is that the original formulation is recovered for $ W_{\mathrm{h}} = V_{\mathrm{h}} $. However, there are many situations, where a different approximation space $ W_{\mathrm{h}} \neq V_{\mathrm{h}} $ may be more natural or even necessary in order to avoid an ill-posed problem.

While the coefficients $ \vec{f} $ are defined through \eqref{eqn:egfem_coeffs}, the approximation space $ W_{\mathrm{h}} $ still needs to be fixed. In particular, the actual setting, e.g., the nonlinear function $ f $, the trial space $ V_{\mathrm{h}} $ and other numerical considerations, influences our choice for the local finite element space. The plausible scenarios that influence how to approximate $ f(u_{\mathrm{h}},\nabla u_{\mathrm{h}}) $ are separated into the following three cases.

\subsubsection*{Case~1 (Exact Reformulation)} 

Let $ K $ denote an element in $ \Omega_{\mathrm{h}} $. At this point, a trial space $ V_{\mathrm{h}} $ based on finite elements has been selected. When we consider the restriction of $ f(u_{\mathrm{h}},\nabla u_{\mathrm{h}}) $ to the element $ K $, it follows that the behavior of $ (f(u_{\mathrm{h}},\nabla u_{\mathrm{h}}))\vert_{K} = f(u_{\mathrm{h}} \vert_{K},\nabla u_{\mathrm{h}} \vert_{K}) $ depends on the choice of $ V_{\mathrm{h}} $. This case assumes that the local restriction can be described exactly by some finite element. This leads us directly to an appropriate approximation space $ W_{\mathrm{h}} $.
\begin{example} \label{eg:exact_egfem}
    Let $ K $ be an element in $ \Omega_{\mathrm{h}} $ and $ \pp_{k}(K) $ denote the space of polynomial functions of at most degree $ k $ on $ K $. Consider $ u_{\mathrm{h}} \vert_{K} \in \pp_{k} $, i.e., $ V_{\mathrm{h}} = \{v : v|_{K} \in \pp_{k} ~\forall K \in \Omega_{\mathrm{h}}\} $. All monomials $ f_{\ell}(u) \mdef u^{\ell} $ with $ \ell \in \nn_{0} $ fulfill $ (f_{\ell}(u_{\mathrm{h}})) \vert_{K} \in \pp_{\ell \cdot k} $, which leads to the choice $ W_{\mathrm{h}} = \{w : w|_{K} \in \pp_{\ell \cdot k} ~\forall K \in \Omega_{\mathrm{h}}\}$. For an example that highlights a major weakness of the \emph{original} group finite element method, consider the nonlinear function defined through $ f(u,\nabla u) \mdef \Vert \nabla u \Vert^{2} $. The local restriction to $ K $ can be characterized through $ f(u_{\mathrm{h}},\nabla u_{\mathrm{h}}) \vert_{K} \in \pp_{2 \, (k-1)} $. Once again, we can choose $ W_{\mathrm{h}} = \{w : w|_{K} \in \pp_{2 \, (k-1)} ~\forall K \in \Omega_{\mathrm{h}}\} $, which exactly describes the behavior of $ f(u,\nabla u) $ on each element for the fixed choice of $ V_{\mathrm{h}} $.
\end{example}%

\subsubsection*{Case~2 (Approximation)} 

Let the trial space $ V_{\mathrm{h}} $ and the approximation space $ W_{\mathrm{h}} $ be fixed. Then, the nonlinear function $ f(u,\nabla u) $ is approximated locally by the finite elements used in the construction of $ W_{\mathrm{h}} $. This case depends on a variety of factors. For example, if Case~1 is not applicable, then it may be acceptable to approximate the nonlinearity directly. Additionally, some problems require other approximations, e.g., errors that arise from evaluating the integration in \eqref{eqn:nonlinear_terms} inexactly. The current case allows for a direct approximation of $ f(u,\nabla u) $ through a given $ W_{\mathrm{h}} $, which takes these situations into consideration.
\begin{example}
    The group finite element method generally falls into this case. Here, the approximation space is the same as the trial space, i.e., $ V_{\mathrm{h}} = W_{\mathrm{h}} $.
\end{example} 

\subsubsection*{Case~3 (Quadrature)} 

For a fixed trial space $ V_{\mathrm{h}} $, we consider the situation, where an exact reformulation of the \emph{discrete} problem is desired but not feasible with existing finite element spaces. In this case, we focus on reconstructing the quadrature rule used in the original formulation. In \mbox{\cite{Wang2015}}, an approach is used for nonlinear terms that embeds the assembly process into a matrix-vector product, where the local quadrature information is stored in the matrix and the vector consists of the nonlinear function evaluated at all quadrature points. For our work, we build upon this idea for use in our extended group finite element formulation. In particular, we construct an approximation space $ W_{\mathrm{h}} $ by defining a set of \emph{discontinuous} basis functions $ \{\eta_{k}^{\mathrm{f}}\}_{k=1}^{N_{\mathrm{f}}} $ that span it and contain the necessary quadrature information. Consider, for example, a $ N_{\mathrm{q}} $-point quadrature rule that approximates the integral
\begin{gather*}
    \int_{K} f(u_{\mathrm{h}}(x),\nabla u_{\mathrm{h}}(x)) \, \phi_{i}(x) \dx{x} \approx \sum_{\ell=1}^{N_{\mathrm{q}}} w_{\ell}^{\mathrm{K}} \, f(u_{\mathrm{h}}(x_{\ell}^{\mathrm{K}}),\nabla u_{\mathrm{h}}(x_{\ell}^{\mathrm{K}})) \, \phi_{i}(x_{\ell}^{\mathrm{K}}),
\end{gather*}
where $ w_{\ell}^{\mathrm{K}} $ denotes the weights of the quadrature rule and $ x_{\ell}^{\mathrm{K}} $ denotes the quadrature nodes located on the element $ K \in \Omega_{\mathrm{h}} $. In order to define a finite element, we need to fix the degrees of freedom and basis functions on an element. Here, the degrees of freedom are given by the quadrature nodes $ \{x_{\ell}^{\mathrm{K}}\}_{\ell=1}^{N_{\mathrm{q}}} $. The basis functions use discrete Dirac-delta functions $\delta_{x_{\ell}^{\mathrm{K}}}^{\mathrm{h}}$, see \mbox{\cite{Scott1973,Houston2012}} for details on the construction of these approximations. By definition, the discrete $ \delta $-functions fulfill the following property:
\begin{gather*}
	\int_{\Omega_{\mathrm{h}}} v(x) \, \delta_{x_{\ell}^{\mathrm{K}}}^{\mathrm{h}}(x) \dx{x} = v(x_{\ell}^{\mathrm{K}})
\end{gather*}
for all $ v \in V_{\mathrm{h}} $. Thus, the basis functions are defined through
\begin{gather} \label{eqn:quadrature_basis}
    \eta_{\ell}^{\mathrm{K}}(x) \mdef w_{\ell}^{\mathrm{K}} \, \delta_{x_{\ell}^{\mathrm{K}}}^{\mathrm{h}}(x).
\end{gather}
Using these basis functions, the nonlinear function $ f(u_{\mathrm{h}},\nabla u_{\mathrm{h}}) $ is approximated by
\begin{gather*}
    (f(u_{\mathrm{h}}(x),\nabla u_{\mathrm{h}}(x)))\vert_{K} \approx (f_{\mathrm{h}}(x))\vert_{K} = \sum_{\ell=1}^{N_{\mathrm{q}}} f_{\ell} \, \eta_{\ell}^{\mathrm{K}}(x),
\end{gather*}
where $ f_{\ell} $ is defined through \eqref{eqn:egfem_coeffs} at the quadrature point $x_{\ell}^{\mathrm{K}}$. Introducing this approximation into the above integral leads to
\begin{gather*}
     \int_{K} f_{\mathrm{h}}(x) \, \phi_{i}(x) \dx{x} = \sum_{\ell=1}^{N_{\mathrm{q}}} f_{\ell} \int_{K} \eta_{\ell}^{\mathrm{K}}(x) \, \phi_{i}(x) \dx{x} = \sum_{\ell=1}^{N_{\mathrm{q}}} f(u_{h}(x_{\ell}^{\mathrm{K}}),\nabla u_{\mathrm{h}}(x_{\ell}^{\mathrm{K}})) \, w_{\ell}^{\mathrm{K}} \, \phi_{i}(x_{\ell}^{\mathrm{K}}),
\end{gather*}
due to the fact that $ \int_{K} w_{\ell}^{\mathrm{K}} \, \delta_{x_{\ell}^{\mathrm{K}}}^{\mathrm{h}}(x) \, \phi_{i}(x) \dx{x} = w_{\ell}^{\mathrm{K}} \, \phi_{i}(x_{\ell}^{\mathrm{K}}) $. Therefore, the quadrature is embedded into the approximation space $ W_{\mathrm{h}} $ and the discrete version of the original problem is conserved. Note that $ W_{\mathrm{h}} $ is constructed by combining the quadrature points $ \{x_{\ell}^{\mathrm{K}}\}_{\ell=1}^{N_{\mathrm{q}}} $ and basis functions $ \{\eta_{\ell}^{\mathrm{K}}\}_{\ell=1}^{N_{\mathrm{q}}} $ for each element $ K \in \Omega_{\mathrm{h}} $. We denote by $ \ii_{k} $ the space described through the basis functions $\{\eta_{\ell}^{\mathrm{K}}\}_{\ell=1}^{N_{\mathrm{q}}}$, where $ k $ identifies the degree of exactness of the underlying quadrature rule. Depending on the chosen quadrature rule, this quickly leads to a large number of unknowns. For example, the quadrature points of symmetrical Gaussian quadrature rules on triangles are located in the interior of the triangle \cite{Dunavant1985}. Therefore, the approximation space constructed using such a quadrature rule consists of $ N_{\mathrm{q}} \cdot N_{\mathrm{el}} $ basis functions and degrees of freedom, where $ N_{\mathrm{q}} $ and $ N_{\mathrm{el}} $ denote the number of quadrature points and the number of elements in the mesh $ \Omega_{\mathrm{h}} $, respectively.
\begin{remark} 
    For this work, we consider any quadrature rule that is defined on the elements of $ \Omega_{h} $. However, it may be relevant -- in the sense of hyper-reduction -- to construct special quadrature rules via empirical cubature on the entire domain that approximate the integrals in the weak formulation, see, for example, \cite{Hernandez2017} for details on such an approach in the context of model order reduction. In this case, an approximation of the original discrete problem is attained instead of an exact reformulation.
\end{remark}%
\begin{example}
    Consider a triangular element $ K $ from $ \Omega_{\mathrm{h}} $. The following quadrature rule is exact for linear functions:
    \begin{gather*}
        \int_{K} f(x) \dx{x} \approx |K| \, f(c),
    \end{gather*}
    where $ c $ denotes the center point of the triangle and $ |K| $ the triangle's area. The associated quadrature finite element consists of a single degree of freedom at the center $ x_{1}^{\mathrm{K}} = c $ with the constructed basis function $ \eta_{1}^{\mathrm{K}} $ defined in \eqref{eqn:quadrature_basis} with $ w_{1}^{\mathrm{K}} = |K| $. The resulting approximation space $ W_{\mathrm{h}} $ is defined on the collection of all the quadrature points over each element and the respective basis functions, i.e., $ W_{\mathrm{h}} = \mathop{\mathrm{span}} \{\eta_{1}^{\mathrm{K}}\}_{K \in \Omega_{\mathrm{h}}} $. Compared to the discontinuous constant finite element, the location of the degree of freedom coincides, while the basis functions differ.
\end{example}

\subsection{Algebraic Equations} \label{subsec:algebraic}

After introducing the approximation space $ W_{\mathrm{h}} $ for the nonlinear function $ f(u_{\mathrm{h}},\nabla u_{\mathrm{h}}) $ in Section \ref{subsec:egfem_choice}, the coefficients $ \vec{f} $ require the interpolated values and/or gradients of the approximate solution $ u_{\mathrm{h}} $ at every degree of freedom belonging to $ \vec{f} $, i.e., $ f_{k} = f(u_{\mathrm{h}}(x_{k}^{\mathrm{f}}),\nabla u_{\mathrm{h}}(x_{k}^{\mathrm{f}})) $. Instead of evaluating each constraint individually, the evaluation can be vectorized with the help of two interpolation operators. First, the evaluation of the approximate solution $ u_{\mathrm{h}} $ on each degree of freedom $ x_{k}^{\mathrm{f}} $ associated to $ \vec{f} $ can be performed at once using 
\begin{gather*}
	(\Pi_{\mathrm{u}}^{\mathrm{f}})_{ij} \mdef \phi_{j}(x_{i}^{\mathrm{f}}), \quad 1 \leq i \leq N_{\mathrm{f}}, 1 \leq j \leq N_{\mathrm{u}}, 
\shortintertext{so that}
    \Pi_{\mathrm{u}}^{\mathrm{f}} \, \vec{u} = \begin{pmatrix} \phi_{1}(x_{1}^{\mathrm{f}}) & \cdots & \phi_{N_{\mathrm{u}}}(x_{1}^{\mathrm{f}}) \\ \vdots & & \vdots \\ \phi_{1}(x_{N_{\mathrm{f}}}^{\mathrm{f}}) & \cdots & \phi_{N_{\mathrm{u}}}(x_{N_{\mathrm{f}}}^{\mathrm{f}}) \end{pmatrix} \, \begin{pmatrix} u_{1} \\ \vdots \\ u_{N_{\mathrm{u}}} \end{pmatrix} = \begin{pmatrix} u_{\mathrm{h}}(x_{1}^{\mathrm{f}}) \\ \vdots \\ u_{\mathrm{h}}(x_{N_{\mathrm{f}}}^{\mathrm{f}}) \end{pmatrix}.
\end{gather*}
Similarly, a third-order tensor interpolates the gradient of the approximate solution $ \nabla u_{\mathrm{h}} $ onto the degrees of freedom through 
\begin{gather*}
	(\Pi_{\nabla \mathrm{u}}^{\mathrm{f}})_{ijk} \mdef (\nabla \phi_{j}(x_{k}^{\mathrm{f}}))_{i}, \quad 1 \leq i \leq n, 1 \leq j \leq N_{\mathrm{u}}, 1 \leq k \leq N_{\mathrm{f}} , 
\shortintertext{so that}
	\Pi_{\nabla \mathrm{u}}^{\mathrm{f}} \cdot_{2} \vec{u} \mdef \sum_{j=1}^{N_{\mathrm{u}}} (\Pi_{\nabla \mathrm{u}}^{\mathrm{f}})_{ijk} \, u_{j} = \begin{pmatrix} | & & | \\ \nabla u_{\mathrm{h}}(x_{1}^{\mathrm{f}}) & \cdots & \nabla u_{\mathrm{h}}(x_{N_{\mathrm{f}}}^{\mathrm{f}}) \\ | & & | \end{pmatrix}.
\end{gather*}
This leads to the ``vectorization'' of the algebraic equations via
\begin{gather*}
    \vec{f} = f(\Pi_{\mathrm{u}}^{\mathrm{f}} \, \vec{u}, \Pi_{\nabla \mathrm{u}}^{\mathrm{f}} \cdot_{2} \vec{u}),
\end{gather*}
which allows for the efficient evaluation of these terms in the numerical implementation.

\subsection{Multi-linear Forms} \label{subsec:egfem_tensors}

The idea of introducing variables as placeholders for nonlinear terms in order to introduce structure to a given problem can be found in many different applications. For example, the idea of quadratization \cite{Gu2011,Kramer2019}, which reformulates nonlinear problems as quadratic-bilinear systems, can be found in the field of model order reduction. Alternatively, ``lifting'' approaches for Newton's method bring structure to nonlinear optimization problems, which often converge faster than the nonlifted variants \cite{Albersmeyer2010}. Similarly, the group finite element formulation and our extended variant also introduce a clear structure to the discrete nonlinear problem \eqref{eqn:discrete_problem}. Here, the nonlinear terms in \eqref{eqn:nonlinear_terms} simplify into multi-linear forms. Note that depending on the nonlinear functions, different forms may arise but can be handled similarly. Our focus in this subsection lies in the structure introduced by interpolating the nonlinearities onto approximation spaces. In particular, nonlinear forms are replaced with multi-linear forms. These simpler forms are represented by appropriate tensors, which prove to be numerically much more efficient in comparison to the assembly of the nonlinear terms.
\begin{notation}
    Let $ \tens{T} $ denote a third-order tensor, $ \mat{A} $ a second-order tensor (matrix) and $ \vec{v} $, $ \vec{w} $ first-order tensors (vectors). We maintain this notation for the remainder of this work. Then, the tensor (dyadic) product is defined through
    \begin{gather*}
        (\vec{v} \otimes \vec{w})_{ij} = v_{i} \, w_{j}.
    \end{gather*}
    The first and second contractions, denoted by $ \cdot $ and $ : $, are defined through
    \begin{gather*}
        (\mat{A} \cdot \vec{v})_{i} = \sum_{j} \mat{A}_{ij} \, v_{j}, \quad        
        (\tens{T} \cdot \vec{v})_{ij} = \sum_{k} \tens{T}_{ijk} \, v_{k}, \quad
        (\tens{T} : \mat{A})_{i} = \sum_{j,k} \tens{T}_{ijk} \, \mat{A}_{jk}.
    \end{gather*}
    Additionally, the first contraction can also be defined along a different dimension, which we denote by appending an index to the $ \cdot $, e.g.,
    \begin{gather*}
    	(\tens{T} \cdot_{2} \vec{v})_{ik} = \sum_{j} \tens{T}_{ijk} \, v_{j}
    \end{gather*}
    as already introduced in Section \ref{subsec:algebraic}. Finally, we denote the Hadamard product, which refers to an element-wise product, through
    \begin{gather*}
    	(\vec{v} \odot \vec{w})_{i} = v_{i} \, w_{i}
    \end{gather*}
    for vectors of the same dimension. From these definitions, it is clear that the tensor product results in a second-order tensor, while the contractions reduce the order of the tensor by one and two orders, respectively. The Hadamard product does not change the order/dimension at all.
\end{notation}
In order to facilitate the discussion regarding the tensor structure of the extended group finite element formulation, we focus on the discrete problem \eqref{eqn:discrete_problem}. For the standard Galerkin approach, which requires the assembly of the nonlinear terms in every iteration of the iterative method, we have the terms $ \vec{K}(a,b,\vec{u}) $, $ \vec{M}(c,\vec{u}) $ and $ \vec{g}(\vec{u}) $ in \eqref{eqn:nonlinear_terms}, where $ \vec{K} $ and $ \vec{M} $ represent generalizations of the classical stiffness and mass matrices. In contrast to the standard formulation, our extended group finite element formulation has a well-defined structure that can be precomputed. Applying the discretization discussed in Section \ref{subsec:egfem_choice} to \eqref{eqn:nonlinear_terms} results in the following discrete forms:
\begin{gather} \label{eqn:egfem_terms}
\begin{aligned}
    \left(\tens{K}^{\mathrm{b}}_{\mathrm{a}} : (\vec{b} \otimes \vec{a})\right)_{i} &\mdef \sum_{k=1}^{N_{\mathrm{a}}} a_{k} \sum_{j=1}^{N_{\mathrm{b}}} b_{j} \left(\int_{\Omega_{\mathrm{h}}} \eta_{k}^{\mathrm{a}} \, \nabla \eta_{j}^{\mathrm{b}} \cdot \nabla \phi_{i} \dx{x}\right), \\
    \left(\mat{M}^{\mathrm{c}} \, \vec{c}\right)_{i} &\mdef \sum_{j=1}^{N_{\mathrm{c}}} c_{j} \left(\int_{\Omega_{\mathrm{h}}} \eta_{j}^{\mathrm{c}} \, \phi_{i} \dx{x}\right), \\
    \left(\mat{G}^{\mathrm{g}} \, \vec{g}\right)_{i} &\mdef \sum_{j=1}^{N_{\mathrm{g}}} g_{j} \left(\int_{\Gamma_{\mathrm{N},\mathrm{h}}} \eta_{j}^{\mathrm{g}} \, \phi_{i} \dx{x}\right)
\end{aligned}
\end{gather}%
with $ (\tens{K}^{\mathrm{b}}_{\mathrm{a}})_{ijk} \mdef \int_{\Omega_{\mathrm{h}}} \eta_{k}^{\mathrm{a}} \, \nabla \eta_{j}^{\mathrm{b}} \cdot \nabla \phi_{i} \dx{x} $, $ (\mat{M}^{\mathrm{c}})_{ij} \mdef \int_{\Omega_{\mathrm{h}}} \eta_{j}^{\mathrm{c}} \, \phi_{i} \dx{x} $ and $ (\mat{G}^{\mathrm{g}})_{ij} \mdef \int_{\Gamma_{\mathrm{N},\mathrm{h}}} \eta_{j}^{\mathrm{g}} \, \phi_{i} \dx{x} $, where $ \{\eta_{k}^{\mathrm{f}}\}_{k=1}^{N_{\mathrm{f}}}$ denotes the nodal basis associated to $ f \in \{a,b,c,g\} $. The interesting development in comparison to the GFEM is the tensor structure that results from approximating $ \vec{K}(a,b,\vec{u}) $. By interpolating both nonlinear functions $ a(u,\nabla u) $ and $ b(u) $, we attain a trilinear form, which is represented with the help of a third-order tensor.
\begin{remark}
    In addition to the tensor notation, it is also possible to express the action of a trilinear form in terms of a matrix-vector product by using matricization and the Kronecker product. In particular, we can rewrite the second contraction with the tensor product as
    \begin{gather*}
        \tens{K}^{\mathrm{b}}_{\mathrm{a}} : (\vec{b} \otimes \vec{a}) = \mat{K}^{\mathrm{b}}_{\mathrm{a}} \, (\vec{a} \mathbin{\otimes_{K}} \vec{b}),
    \end{gather*}
    where $ \mat{K}^{\mathrm{b}}_{\mathrm{a}} $ denotes the mode-2 matricization of the tensor $ \tens{K}^{\mathrm{b}}_{\mathrm{a}} $ and $ \mathbin{\otimes_{K}} $ the Kronecker product. However, due to the size of the systems that result from the mode-2 matricization and Kronecker product, the tensor formulation is much more efficient for the purpose of numerical implementation.
\end{remark}

\section{Numerical Methods} \label{sec:numeric}

The advantages of the extended group finite element formulation become even more clear, when the discrete systems, which must be solved in each iteration, are compared to the standard Galerkin approach. In practice, there are two common iterative methods for solving nonlinear systems such as \eqref{eqn:discrete_problem}, see \cite{Logg2012}. The first method uses the Picard iteration and is also known as the method of successive iterations. The idea is to substitute the previous iteration into the nonlinear terms, which delivers a linear system in terms of the new iteration. The second method, Newton's method, is a well-known method for solving nonlinear systems, which exploits gradient information in order to improve the (local) convergence rate. However, this comes at the cost of evaluating the Jacobian matrix, which is also assembled using numerical integration. While Newton's method generally converges more quickly, the Picard iterations are simpler to implement. In light of its simplicity, the benefits of the extended group finite element method are demonstrated using the Picard iterations in Section \ref{sec:results}. The results are qualitative and, therefore, expected to be similar when using Newton's method to solve the nonlinear systems.

The standard Galerkin approach solves the following nonlinear form
\begin{gather} \label{eqn:SGA_root}
    \vec{F}(\vec{u}) = \vec{K}(a,b,\vec{u}) + \vec{M}(c,\vec{u}) - \vec{d} - \vec{h} + \vec{g}(\vec{u}) = \vec{0}
\end{gather}
with the nonlinear terms defined in \eqref{eqn:nonlinear_terms}, while the extended group finite element formulation leads to
\begin{gather} \label{eqn:EGFEM_root}
    \vec{\tilde{F}}(\vec{z}) = \begin{pmatrix} 
        \tens{K}_{\mathrm{a}}^{\mathrm{b}} : (\vec{b} \otimes \vec{a}) + \mat{M}^{\mathrm{c}} \, \vec{c} - \vec{d} - \vec{h} + \mat{G}^{\mathrm{g}} \, \vec{g} \\ 
        \vec{a} - a(\Pi_{\mathrm{u}}^{\mathrm{a}} \, \vec{u}, \Pi_{\nabla \mathrm{u}}^{\mathrm{a}} \cdot_{2} \vec{u}) \\ 
        \vec{b} - b(\Pi_{\mathrm{u}}^{\mathrm{b}} \, \vec{u}) \\ 
        \vec{c} - c(\Pi_{\mathrm{u}}^{\mathrm{c}} \, \vec{u}, \Pi_{\nabla \mathrm{u}}^{\mathrm{c}} \cdot_{2} \vec{u}) \\ 
        \vec{g} - g(\Pi_{\mathrm{u}}^{\mathrm{g}} \,\vec{u}) 
    \end{pmatrix} = \vec{0}
\end{gather}
with $ \vec{z} \mathop{:=} (\vec{u},\vec{a},\vec{b},\vec{c},\vec{g})^{\transpose} $ and the new interpolated terms defined in \eqref{eqn:egfem_terms}. The point-wise definitions of the coefficient vectors $ \vec{a} $, $ \vec{b} $, $ \vec{c} $ and $ \vec{g} $ assume that the associated function evaluations of $ a $, $ b $, $ c $ and $ g $ are vectorized, e.g., $ f_{k} = f(u_{\mathrm{h}}(x_{k}^{\mathrm{f}}),\nabla u_{\mathrm{h}}(x_{k}^{\mathrm{f}})) $, $ 1 \leq k \leq N_{\mathrm{f}} $, which is realized with the help of the interpolation operators $ \Pi_{\mathrm{u}}^{\mathrm{f}} $ and $ \Pi_{\nabla \mathrm{u}}^{\mathrm{f}} $ from Section \ref{subsec:algebraic}, for $ f \in \{a,b,c,g\} $.
\begin{remark}
	Given the general setting of \eqref{eqn:SGA_root} and \eqref{eqn:EGFEM_root}, it is difficult to ensure the existence of solutions. In light of this, there may be constellations of $ V_{\mathrm{h}} $ and $ W_{\mathrm{h}} $ such that the existence of a solution and/or convergence to the solution are lost. Such phenomena occur for mixed problems, e.g., checkerboard instability in the Stokes problem \cite{Ern2004}. Therefore, each concrete problem and discretization must be investigated carefully. For the sake of discussion, we assume that all forms presented in this work are well-defined.
\end{remark}

\subsubsection*{Picard Iteration}

The form of the Picard iterations is problem-dependent. In general, the desired form of the discrete nonlinear problem is given by
\begin{gather*}
    \mat{A}(\vec{u}) \, \vec{u} - \vec{r}(\vec{u}) = \vec{0}.
\end{gather*}
In this setting, the Picard iteration is attained by substituting the previous iteration into the nonlinear terms, i.e.,
\begin{gather*}
    \vec{u}^{n+1} = \mat{A}(\vec{u}^{n})^{-1} \, \vec{r}(\vec{u}^{n}).
\end{gather*}
Because the general problem \eqref{eqn:EGFEM_root} using the EGFEM formulation does not suggest how an iteration of the above form could be attained, let us consider the following example.
\begin{example} \label{eg:picard_iteration}
Simplify \eqref{eqn:nonlinear_pde} by letting $ b(u) = u $. Then, the original term $ \vec{K}(a,b,\vec{u}) $ in \eqref{eqn:SGA_root} can be written as $ \mat{K}(a,\vec{u}) \, \vec{u} $ with
\begin{gather*}
    (\mat{K}(a,\vec{u}))_{ij} \mdef \int_{\Omega_{\mathrm{h}}} a(u_{\mathrm{h}},\nabla u_{\mathrm{h}}) \, \nabla \phi_{j} \cdot \nabla \phi_{i} \dx{x},
\end{gather*}
while the tensor form $ \tens{K}^{\mathrm{b}}_{\mathrm{a}} : (\vec{b} \otimes \vec{a}) $ in \eqref{eqn:EGFEM_root} simplifies to $ \tens{K}_{\mathrm{a}} : (\vec{u} \otimes \vec{a}) = (\tens{K}_{\mathrm{a}} \cdot \vec{a}) \cdot \vec{u} $ with
\begin{gather*}
	(\tens{K}_{\mathrm{a}})_{ijk} \mdef \int_{\Omega_{\mathrm{h}}} \eta_{k}^{\mathrm{a}} \, \nabla \phi_{j} \cdot \nabla \phi_{i} \dx{x}.
\end{gather*}
This leads to the following iteration for the standard Galerkin approach:
\begin{gather*}
    \vec{u}^{n+1} = \mat{K}(a,\vec{u}^{n})^{-1} \left(\vec{d} + \vec{h} - \vec{g}(\vec{u}^{n}) - \vec{M}(c,\vec{u}^{n})\right).
\end{gather*}
In contrast, the (extended) group finite element formulation solves
\begin{align*}
    \vec{u}^{n+1} &= (\tens{K}_{\mathrm{a}} \cdot \vec{a}^{n})^{-1} \left(\vec{d} + \vec{h} - \mat{G}^{\mathrm{g}} \, \vec{g}^{n} - \mat{M}^{\mathrm{c}} \, \vec{c}^{n}\right)
\end{align*}
before updating the nonlinear variables
\begin{align*}
    \vec{f}^{n+1} &= f(\Pi_{\mathrm{u}}^{\mathrm{f}} \, \vec{u}^{n+1}, \Pi_{\nabla \mathrm{u}}^{\mathrm{f}} \cdot_{2} \vec{u}^{n+1}), &
    \vec{g}^{n+1} &= g(\Pi_{\mathrm{u}}^{\mathrm{g}} \, \vec{u}^{n+1}), & f &\in \{a,c\}
\end{align*}
in each iteration.
\end{example}%
Although Example \ref{eg:picard_iteration} considers a special class of problems, it still illustrates that the reformulated systems avoid assembly of the nonlinear terms in each iteration. Instead, the nonlinear functions are evaluated point-wise, which is often computationally less expensive than the numerical integration of the nonlinear forms.

\subsubsection*{Newton Iteration}

The iterations defined by Newton's method can be attained through two different approaches \cite{Logg2012}. The standard approach is to apply Newton's method to the discrete nonlinear problem \eqref{eqn:discrete_problem}, while alternatively it can also be applied to the continuous nonlinear PDE \eqref{eqn:nonlinear_pde}. In the end, both approaches lead to the same system. The iterates are given using the following update:
\begin{gather*}
    \vec{u}^{n+1} = \vec{u}^{n} - \mathrm{D}_{\vec{u}}\vec{F}(\vec{u}^{n})^{-1} \, \vec{F}(\vec{u}^{n}),
\end{gather*}
where $ \mathrm{D}_{\vec{u}}\vec{F} $ denotes the Jacobian matrix. For the standard Galerkin formulation, this system uses \eqref{eqn:SGA_root} with the Jacobian matrix
\begin{gather*}
    \mathrm{D}_{\vec{u}}\vec{F}(\vec{u}) = \mathrm{D}_{\vec{u}}\vec{K}(a,b,\vec{u}) + \mathrm{D}_{\vec{u}}\vec{M}(c,\vec{u}) + \mathrm{D}_{\vec{u}}\vec{g}(\vec{u}).
\end{gather*}
In order to evaluate the Jacobian matrix, assembly using numerical integration is necessary. This leads to a noticeable increase in the computational complexity of each iteration in comparison to the Picard iterations although the Newton iterations are expected to converge more quickly.

The extended group finite element formulation iterates over $ \vec{z} = (\vec{u},\vec{a},\vec{b},\vec{c},\vec{g})^{\transpose} $ using
\begin{gather*}
    \vec{z}^{n+1} = \vec{z}^{n} - \mathrm{D}_{\vec{z}}\vec{\tilde{F}}(\vec{z}^{n})^{-1} \, \vec{\tilde{F}}(\vec{z}^{n}),
\end{gather*}
with $ \vec{\tilde{F}} $ defined in \eqref{eqn:EGFEM_root} and the Jacobian matrix
\begin{gather*}
    \mathrm{D}_{\vec{z}}\vec{\tilde{F}}(\vec{z}) = \begin{pmatrix} \mathbb{0} & \tens{K}_{\mathrm{a}}^{\mathrm{b}} \cdot_{2} \vec{b} & \tens{K}_{\mathrm{a}}^{\mathrm{b}} \cdot \vec{a} & \mat{M}^{\mathrm{c}} & \mat{G}^{\mathrm{g}} \\ \mathrm{D}_{\vec{u}}a(\Pi_{\mathrm{u}}^{\mathrm{a}} \, \vec{u},\Pi_{\nabla \mathrm{u}}^{\mathrm{a}} \cdot_{2} \vec{u}) & \mathbb{1} & \mathbb{0} & \mathbb{0} & \mathbb{0} \\ \mathrm{D}_{\vec{u}}b(\Pi_{\mathrm{u}}^{\mathrm{b}} \, \vec{u}) & \mathbb{0} & \mathbb{1} & \mathbb{0} & \mathbb{0} \\ \mathrm{D}_{\vec{u}}c(\Pi_{\mathrm{u}}^{\mathrm{c}} \, \vec{u},\Pi_{\nabla \mathrm{u}}^{\mathrm{c}} \cdot_{2} \vec{u}) & \mathbb{0} & \mathbb{0} & \mathbb{1} & \mathbb{0} \\ \mathrm{D}_{\vec{u}}g(\Pi_{\mathrm{u}}^{\mathrm{g}} \, \vec{u}) & \mathbb{0} & \mathbb{0} & \mathbb{0} & \mathbb{1} \end{pmatrix}.
\end{gather*}
In contrast to the standard Galerkin formulation, the Jacobian matrix here does not require any additional integration. Instead, the only integrated portions depend on the precomputed forms $ \tens{K}^{\mathrm{b}}_{\mathrm{a}} $, $ \mat{M}^{\mathrm{c}} $ and $ \mat{G}^{\mathrm{g}} $. In addition, the Jacobian matrices $ \mathrm{D}_{\vec{u}} a $, $ \mathrm{D}_{\vec{u}} b $, $ \mathrm{D}_{\vec{u}} c $ and $ \mathrm{D}_{\vec{u}} g $ associated to the nonlinear functions $ a $, $ b $, $ c $ and $ g $ are defined point-wise. 
\begin{remark}
	At this point, we would like to call to attention the precarious balance between the increased system size and the benefits of precomputing the nonlinear forms. This is quite similar to the problem faced by some model order reduction techniques, where new variables have to be introduced before reducing the model \cite{Gu2011,Kramer2019}. The SGA solves a smaller system in comparison to the EGFEM but also requires numerical integration in order to evaluate $ \vec{F} $ and $ \mathrm{D}_{\vec{u}} \vec{F} $. In contrast, the extended group finite element method leads to a much larger system, but the nonlinear form $ \vec{\tilde{F}} $ and its Jacobian matrix $ \mathrm{D}_{\vec{z}}\vec{\tilde{F}} $ can be efficiently evaluated. These aspects must be taken into account, when considering which approach to use.
\end{remark}%
\begin{remark}
    In some cases, it may be necessary to use some kind of dampening or continuation method in order to ensure convergence of the iterations. Because this depends entirely on the nonlinearities and the formulation of the discrete problem, there is no general guideline to determine if such methods are necessary. An example of such a problem is presented in Section \ref{subsec:superconductivity}.
\end{remark}

\section{Results} \label{sec:results}

In this section, we present a comprehensive performance study of the extended group finite element method (EGFEM) for a variety of examples. The examples serve to investigate different aspects of the extensions presented in Section \ref{sec:egfem}, especially in comparison to the standard Galerkin approach (SGA) and original group finite element method (GFEM). In addition to the advantages and disadvantages of the tensor formulation, the exact reformulation of the EGFEM contrasts the approximative nature of the GFEM. Each example begins with a brief description of the nonlinear problem, where the structure is presented with respect to the general form \eqref{eqn:nonlinear_pde}. Afterwards, the discrete systems related to the SGA as well as the original and extended GFEMs are identified. Finally, the important aspects of the given example are discussed and supported with the help of figures and tables. The wide variety of the examples shows the versatility and applicability of our extended group finite element method for nonlinear problems. After presenting all the different examples, we conclude this section with a discussion about the results, including some guidelines to consider when applying the extended group finite element method.

All simulations are performed in Matlab\footnote{All results were computed on a machine with an i7-8700 and 32 GB of RAM using Matlab R2019b.} using Sandia's Tensor Toolbox for the tensor calculations and Gmsh for the mesh generation \cite{TTB_Software,TTB_Sparse,Geu09}. For simplicity, the solution $ u $ is discretized using linear finite elements in all of the examples. The assembly of the system is performed with the help of symmetrical Gaussian quadrature rules for triangles \cite{Dunavant1985}. In light of the additional complexity associated to Newton's method, we use Picard iterations in order to solve the nonlinear problems to a tolerance of $ \mathcal{O}(10^{-12}) $ in the difference between iterations. The solutions of the linear systems that have to be solved in each iteration are computed with the help of Matlab's backslash operator. Also note that the runtimes presented refer only to the online computational times, i.e., the time spent in the iteration method including the assembly of any terms that cannot be precomputed. The time spent computing the system tensors, matrices and vectors as well as the interpolation operators is neglected, since most relevant applications such as many-query and/or real-time optimal control perform these computations in an offline phase. However, it should be noted that these calculations can be time-consuming.

\subsection{Quadratic Nonlinearity} \label{subsec:quadratic}

A simple example of a nonlinear problem is the Poisson problem with an additional quadratic term:
\begin{align*}
    -\Delta u(x) + u^{2}(x) &= d(x), & x &\in \Omega \mdef [0,1] \times [0,1], \\
    u(x) &= u_{\mathrm{D}}(x) = x_{1} \, x_{2} \, (x_{1} + x_{2}), & x &\in \Gamma.
\end{align*}
The source term $ d $ is chosen according to the method of manufactured solutions \cite{Roache2019}, where the desired solution is given through the extension of $ u_{\mathrm{D}} $ onto the entire domain $ \bar{\Omega} $, i.e.,
\begin{gather*}
    d(x) = -2 \, (x_{1} + x_{2}) + x_{1}^{2} \, x_{2}^{2} \, (x_{1} + x_{2})^{2}
\end{gather*}
In relation to \eqref{eqn:nonlinear_pde}, this example is given through
\begin{align*}
    a(u,\nabla u) &= 1, & b(u) &= u, & c(u,\nabla u) &= u^{2}
\end{align*}
with $ \Gamma_{\mathrm{D}} = \Gamma $ and $ \Gamma_{\mathrm{N}} = \emptyset $. After discretizing with linear Lagrange elements, the standard Galerkin approach using Picard iterations delivers the following nonlinear discrete problem:
\begin{gather*}
    \mat{K} \, \vec{u}^{n+1} + \vec{M}(c,\vec{u}^{n}) = \vec{d}
\end{gather*}
with the stiffness matrix $ \mat{K} $. An interesting aspect of this example is that the nonlinearity $ c $ can be handled in two different ways. Firstly, it can be treated directly as a trilinear form, e.g.,
\begin{gather*}
    \int_{\Omega_{\mathrm{h}}} u_{\mathrm{h}}^{2} \, \phi_{i} \dx{x} = \sum_{j=1}^{N_{\mathrm{u}}} u_{j} \, \sum_{k=1}^{N_{\mathrm{u}}} u_{k} \left(\int_{\Omega_{\mathrm{h}}} \phi_{k} \, \phi_{j} \, \phi_{i} \dx{x}\right) \mathop{=:} \tens{M} : (\vec{u} \otimes \vec{u}).
\end{gather*}
This eliminates the need for the repetitive assembly of the nonlinear term in each iteration, i.e., the discrete problem is given through
\begin{gather*}
    \mat{K} \, \vec{u}^{n+1} + \tens{M} : (\vec{u}^{n} \otimes \vec{u}^{n}) = \vec{d}.
\end{gather*}
Note that this approach is only possible for polynomial terms, since these can be formulated as multi-linear forms. However, it is only viable for low-degree polynomials, since the size of the resulting tensor form grows exponentially with respect to the degree of the polynomial. Alternatively, the (extended) group finite element method introduces the variable $ \vec{c} $, which is defined point-wise. This leads to the replacement of the second term $ \vec{M}(c,\vec{u}) $ and introduction of an additional equation
\begin{gather*}
    \mat{K} \, \vec{u}^{n+1} + \mat{M}^{\mathrm{c}} \, \vec{c}^{n} = \vec{d}, \quad 
    \vec{c}^{n+1} = \left(\Pi_{\mathrm{u}}^{\mathrm{c}} \, \vec{u}^{n+1}\right) \odot \left(\Pi_{\mathrm{u}}^{\mathrm{c}} \, \vec{u}^{n+1}\right),
\end{gather*}
where the interpolation operator $ \Pi_{\mathrm{u}}^{\mathrm{c}} $ is equal to the identity in the case of the original group finite element formulation. While the GFEM chooses $ W_{\mathrm{h}} = V_{\mathrm{h}} $, our EGFEM allows for a more general choice. The space of piece-wise quadratic functions ($ \pp_{2} $) exactly describes the nonlinearity. Hence, the nonlinear function $ c $ falls into Case~1 of Section \ref{subsec:egfem_choice} for the extended group finite element method with quadratic Lagrange elements for $ W_{\mathrm{h}} $. For the sake of comparison, we also consider Case~3 by using a quadrature rule that is exact for cubic polynomials ($ \ii_{3} $).

\begin{figure} \centering
	\includegraphics{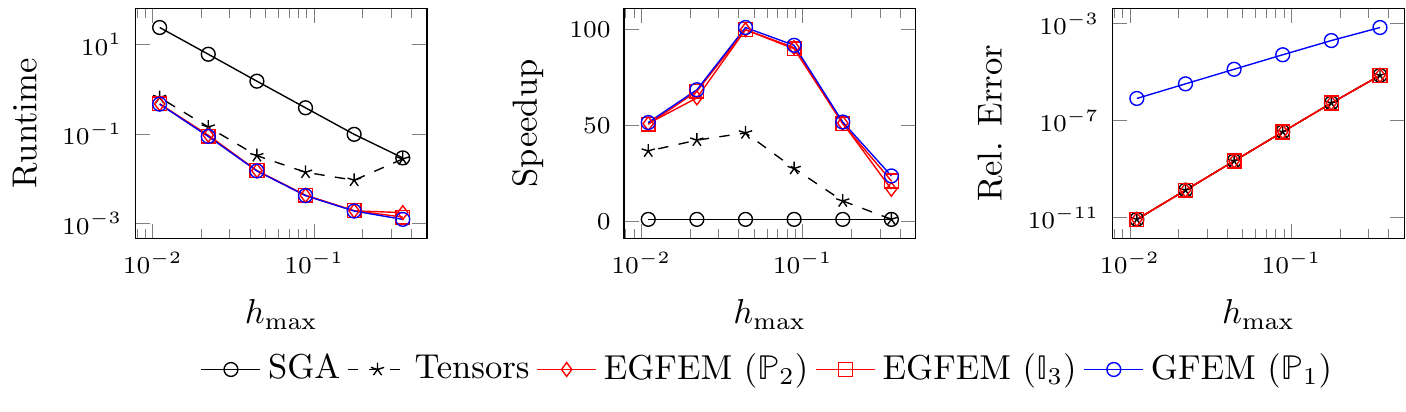}

    \caption{An overview of the performance of each approach for the quadratic example using a regular mesh that is successively refined through splitting. \textsl{From left to right:} The amount of time to solve the system for increasingly finer meshes, the speedup describes the ratio between the reference runtime (SGA) and the considered approach, and the relative error with respect to the approximate $ L^{2} $-norm from the exact solution.}
    \label{fig:quadratic}
\end{figure}

\begin{figure} \centering
	\includegraphics{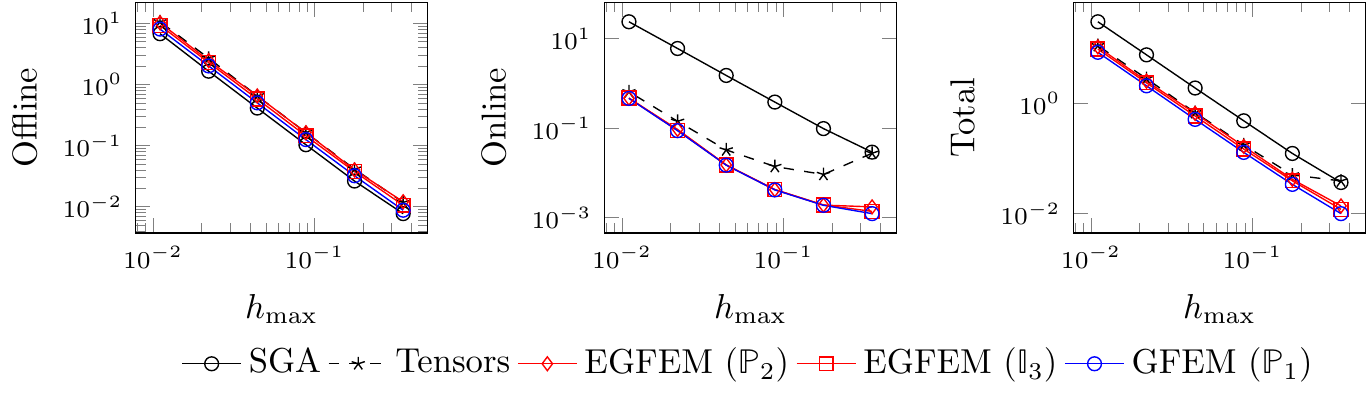}

    \caption{An overview of the time spent computing the solution for each method. \textsl{From left to right:} The offline time includes the assembly of any terms that can be precomputed, e.g., the stiffness matrix and the source as well as the third-order tensor and interpolation operator for the EGFEM, the online time is the time spent solving the iterations to the given tolerance, and the total time combines both the offline and online times together.}
    \label{fig:quadratic_runtime}
\end{figure}

A comprehensive overview of the numerical performance of each approach can be found in Figure \ref{fig:quadratic}. When comparing the runtime and, therefore, also the speedup, we see that the tensor formulation becomes more efficient as the size of the problem grows. However, the original and extended group finite element formulations immediately show an improvement that delivers a speedup between one to two orders of magnitude. In Figure \ref{fig:quadratic_runtime}, we also compare the offline time and total time spent solving the discrete problems to the online runtime. The offline computations take longer for the reformulations, which is to be expected. When the combined time is considered, we see that the gap between the standard approach and the original and extended group finite element methods shrinks. However, there is still an advantage, which grows as the number for solves increases. The important difference between the GFEM and the EGFEM can be seen in the relative error with respect to the exact solution. Here, the GFEM converges slower and with a larger error than the extended formulation using $ \pp_{2} $ and $ \ii_{3} $. This lies on the approximative nature of the group finite element formulation. In particular, the GFEM interpolates the nonlinearity onto the space of piece-wise linear functions ($ \pp_{1} $). By choosing the space of piece-wise quadratic functions ($ \pp_{2} $) or the cubic quadrature rule ($ \ii_{3} $), the EGFEM preserves the nature of the nonlinearity. Therefore, the extended group finite element formulation in these cases are numerically equivalent to the standard Galerkin approach, while exhibiting the same increase in performance as the original formulation. 

\begin{figure} \centering
	\includegraphics{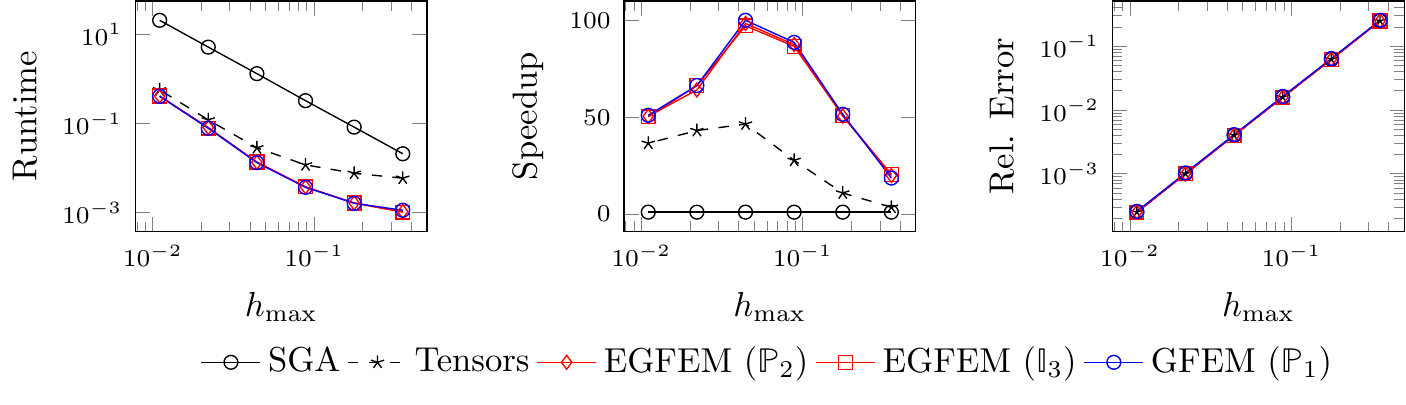}

    \caption{An overview of the performance of each approach for the quadratic example using a different choice for the Dirichlet boundary condition. \textsl{From left to right:} The amount of time to solve the system for increasingly finer meshes, the speedup describes the ratio between the reference runtime (SGA) and the considered approach, and the relative error with respect to the approximate $ L^{2} $-norm from the exact solution.}
    \label{fig:quadraticData}
\end{figure}

In order to better understand the interpolation error, consider a different choice for the Dirichlet boundary condition $ u_{\mathrm{D}} $ and source term $ d $. For example, the same computations as before now display a different behavior in the error convergence for 
\begin{gather*}
	u_{\mathrm{D}}(x) = \sin(2 \, \pi \, x_{1}) \, \sin(2 \, \pi \, x_{2})
\end{gather*}
with the source $ d(x) = -\Delta u_{\mathrm{D}}(x) + u_{\mathrm{D}}(x)^{2} $. The results are very interesting, since there appears to be no perceivable difference between the various approximations of the nonlinearity (cf. Figure \ref{fig:quadraticData}). In our opinion, the approximation error based on the choice of $ V_{\mathrm{h}} $ and/or the additional errors introduced in inexactly integrating the source term are -- in this case -- much more significant than the error based on the approximation of the nonlinear function $ c(u,\nabla u) $. This phenomenon reoccurs in other examples considered in this section and, therefore, plays a vital role, when choosing the approximation space $ W_{\mathrm{h}} $.

\subsection{Burgers' Equation} \label{subsec:burgers}

A common example used in many literature sources for illustrating the group finite element formulation is Burgers' equation, see, for example, \cite{Christie1981,Fletcher1983,Dickinson2010}. Here, we consider the viscous two-dimensional Burgers' equation with homogeneous Dirichlet boundary values on the time interval $ [0,T] $ with $ T > 0 $, i.e., 
\begin{align*}
    \partial_{t} u - \nu \, \Delta u + u \, \partial_{x_{1}} u + u \, \partial_{x_{2}} u &= d, & x &\in \Omega \mdef [0,1] \times [0,1], \; t \in [0,T]
\end{align*}
with the initial value $ u(x,0) = u_{0}(x) $ for all $ x \in \Omega $ and $ u(x,t) = 0 $ on the boundary. Adopting the same setting as \cite{Dickinson2010}, the initial value $ u_{0} $ and the source $ d $ are chosen so that the exact solution is given through the following function:
\begin{gather*}
    u(x,t) = 10 \, x_{1} \, x_{2} \, (x_{1} - 1) \, (x_{2} - 1) \left[\sin(2 \, x_{1} \, t) \, e^{-\frac{1}{2} \, t} + \cos(x_{2} \, t) \, e^{-\frac{1}{4} \, t} + \sin(x_{1} \, x_{2} \, t) \, e^{-t}\right],
\end{gather*}
while we choose $ \nu = 1 $ and $ T = 1 $. By taking advantage of the fact that the nonlinear terms can be rewritten as 
\begin{gather*}
    u \, \partial_{x_{1}} u + u \, \partial_{x_{2}} u = \frac{1}{2} \left(\partial_{x_{1}} (u^{2}) + \partial_{x_{2}} (u^{2})\right),
\end{gather*}
we apply integration by parts in order to attain the following semi-discrete form for the standard Galerkin approach:
\begin{gather*}
    \mat{M} \, \partial_{t} \vec{u}(t) + \nu \, \mat{K} \, \vec{u}(t) - \frac{1}{2} \, \vec{N}(\vec{u}(t)) = \vec{d}(t),
\end{gather*}
with
\begin{gather*}
    \vec{N}(\vec{u}(t))_{i} \mdef \int_{\Omega_{\mathrm{h}}} u_{h}(x,t)^{2} \left(\partial_{x_{1}} \phi_{i}(x) + \partial_{x_{2}} \phi_{i}(x)\right) \dx{x}.
\end{gather*}
However, since the action $ \vec{N}(\vec{u}(t)) $ can also be described with the help of a multi-linear form, we also introduce the following tensor-based formulation, which does not require repeated assembly of the system:
\begin{gather*}
    \mat{M} \, \partial_{t} \vec{u}(t) + \nu \, \mat{K} \, \vec{u}(t) - \frac{1}{2} \, \tens{N} : (\vec{u}(t) \otimes \vec{u}(t)) = \vec{d}(t),
\end{gather*}
where the tensor $ \tens{N} $ is defined through
\begin{gather*}
    \tens{N}_{ijk} \mdef \int_{\Omega_{\mathrm{h}}} \phi_{k}(x) \, \phi_{j}(x) \left( \partial_{x_{1}} \phi_{i}(x) + \partial_{x_{2}} \phi_{i}(x) \right) \dx{x}.
\end{gather*}
The extended and original group finite element methods interpolate the term $ f(u) = u^{2} $ onto the approximation space $ W_{\mathrm{h}} $ and $ V_{\mathrm{h}} $, respectively. This leads to the semi-discrete form
\begin{gather*}
    \mat{M} \, \partial_{t} \vec{u}(t) + \nu \, \mat{K} \, \vec{u}(t) - \frac{1}{2} \, \mat{N}^{\mathrm{f}} \, \vec{f}(t) = \vec{d}(t), \quad \vec{f}(t) = (\Pi_{\mathrm{u}}^{\mathrm{f}} \vec{u}) \odot (\Pi_{u}^{f} \vec{u})
\end{gather*}
with
\begin{gather*}
    \mat{N}^{\mathrm{f}}_{ij} \mdef \int_{\Omega_{h}} \eta_{j}^{\mathrm{f}}(x) \, \left( \partial_{x_{1}} \phi_{i}(x) + \partial_{x_{2}} \phi_{i}(x) \right) \dx{x}.
\end{gather*}
We apply a semi-implicit discretization for the temporal component in order to attain an iteration rule similar to the Picard iterations, i.e., by using the previous iterate for the nonlinear terms and the new iterate for all other terms with the time step $ \delta t $. This leads to 
\begin{align*}
	\mat{M} \, (\vec{u}^{n+1} - \vec{u}^{n}) + \delta t \left(\nu \, \mat{K} \, \vec{u}^{n+1} - \frac{1}{2} \, \vec{N}(\vec{u}^{n}) - \vec{d}^{n+1}\right) = 0
\end{align*}
for the standard Galerkin approach. By using a tensor representation for the multi-linear form, the SGA can also be defined through
\begin{align*}
	\mat{M} \, (\vec{u}^{n+1} - \vec{u}^{n}) + \delta t \left(\nu \, \mat{K} \vec{u}^{n+1} - \frac{1}{2} \, \tens{N} : (\vec{u}^{n} \otimes \vec{u}^{n}) - \vec{d}^{n+1}\right) = 0.
\end{align*}
Finally, the (extended) group finite element formulation solves the following system in each time step:
\begin{gather*}
	\mat{M} \, (\vec{u}^{n+1} - \vec{u}^{n}) + \delta t \left(\nu \, \mat{K} \vec{u}^{n+1} - \frac{1}{2} \, \mat{N}^{\mathrm{f}} \, \vec{f}^{n} - \vec{d}^{n+1}\right) = 0, \quad \vec{f}^{n+1} = (\Pi_{\mathrm{u}}^{\mathrm{f}} \vec{u}^{n+1}) \odot (\Pi_{\mathrm{u}}^{\mathrm{f}} \vec{u}^{n+1}).
\end{gather*}
The equation for $ \vec{u}^{n+1} $ is linear, since the nonlinear terms are evaluated using the previous iteration. Afterwards, the nonlinear coefficient $ \vec{f}^{n+1} $ can be updated with the help of the new iterate $ \vec{u}^{n+1} $.

\begin{figure} \centering
	\includegraphics{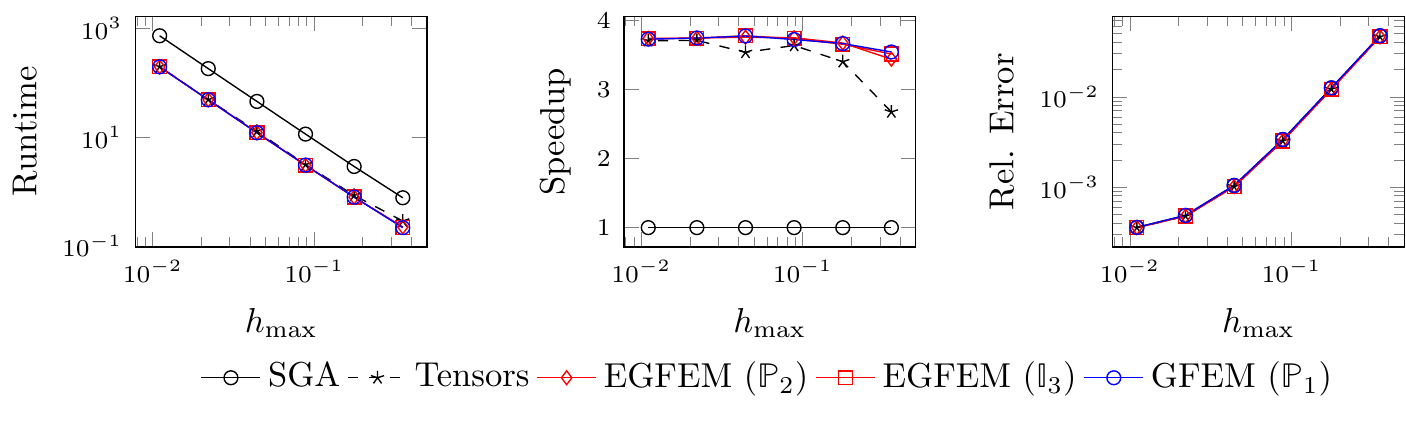}

    \caption{An overview of the performance of each approach for the viscous Burgers' equation using a regular mesh that is successively refined through splitting. \textsl{From left to right:} The average amount of time to solve the problem for increasingly finer meshes, the speedup describes the ratio between the reference runtime (SGA) and the considered approach, and the relative error with respect to the approximate $ L^{2} $-norm from the exact solution.}
    \label{fig:burgers}
\end{figure}

Figure \ref{fig:burgers} compares the performance of each method, when using an equidistant time mesh with $ \delta t = 10^{-2} $ for the time discretization. We see that the tensor-based formulation improves the performance, since the computational overhead associated to assembling the system in each iteration is removed. Additional performance gains are achieved by introducing the (extended) group finite element formulation, which relies on matrices instead of a third-order tensor for the nonlinear terms. This difference leads to the slight improvements between the (extended) group finite element approach and the tensor formulation. In contrast to the first benchmark problem, we see that the difference in accuracy between the original and extended group finite element formulations is numerically negligible. While the GFEM interpolates the nonlinearity onto a piece-wise linear approximation space ($ \pp_{1} $), the EGFEM exactly describes the nonlinearity in the discrete setting by using the finite element space composed of piece-wise quadratic functions ($ \pp_{2} $) or cubic quadrature rule ($ \ii_{3} $). When compared to the exact solution, the error differences are minimal. Based on these and following results, we see that Case~2 in Section \ref{subsec:egfem_choice} is suitable in many different situations. However, as we have also seen in Section \ref{subsec:quadratic}, the EGFEM can lead to noticeably more accurate results with similar performance gains. In the end, the choice of the approximation space depends on the setting and problem to be solved. 

\subsection{Superconductivity} \label{subsec:superconductivity}

For this example, we consider a simplified Ginzburg-Landau model for superconductivity as presented in \cite{Hintermueller2005,Ito1996}. The semi-linear model is given by
\begin{align*}
    -\nu \, \Delta u(x) + u^{3}(x) + u(x) &= d(x), & x &\in \Omega \mdef [0,1] \times [0,1], \\
    u(x) &= u_{\mathrm{D}}(x), & x &\in \Gamma.
\end{align*}
We use the method of manufactured solutions in order to determine the source $ d $ as
\begin{gather*}
    d(x) = -\nu \, \Delta u_{\mathrm{D}}(x) + u_{\mathrm{D}}(x)^{3} + u_{\mathrm{D}}(x).
\end{gather*}
The Dirichlet boundary condition is chosen as
\begin{gather*}
 u_{\mathrm{D}}(x) = \frac{1}{6} \, \sin(2 \, \pi \, x_{1}) \, \sin(2 \, \pi \, x_{2}) \, \exp(2 \, x_{1}).
\end{gather*}
In the framework of \eqref{eqn:nonlinear_pde}, we have 
\begin{align*}
    a(u,\nabla u) &= \nu, & b(u) &= u, & c(u,\nabla u) &= u^{3} + u
\end{align*} 
with $ \Gamma_{\mathrm{D}} = \Gamma $ and $ \Gamma_{\mathrm{N}} = \emptyset $. This example allows us to investigate two different aspects, which can influence the convergence of the iterations. On the one side, the parameter $ \nu $ controls the importance of the nonlinearity $ c $ for $ 0 < \nu \ll 1 $. On the other hand, the nonlinearity $ c $ can be decomposed in the following three different but analytically equivalent ways:
\begin{subequations}
\begin{align}
    \int_{\Omega_{\mathrm{h}}} c(u_{\mathrm{h}},\nabla u_{\mathrm{h}}) \, \phi_{i} \dx{x} &= \mat{M}(\tilde{c},\vec{u}) \, \vec{u}, & \tilde{c}(u,\nabla u) &= u^{2} + 1, \\
     \int_{\Omega_{\mathrm{h}}} c(u_{\mathrm{h}},\nabla u_{\mathrm{h}}) \, \phi_{i} \dx{x} &= \mat{M} \, \vec{u} + \vec{M}(\tilde{c},\vec{u}), & \tilde{c}(u,\nabla u) &= u^{3}, \\
     \int_{\Omega_{\mathrm{h}}} c(u_{\mathrm{h}},\nabla u_{\mathrm{h}}) \, \phi_{i} \dx{x} &= \vec{M}(c,\vec{u}),
\end{align}
\end{subequations}
where the (weighted) mass matrix $ \mat{M}(\tilde{c},\vec{u}) $ is defined through
\begin{gather*}
    \mat{M}(\tilde{c},\vec{u}) \mdef \int_{\Omega_{\mathrm{h}}} \tilde{c}(u_{\mathrm{h}},\nabla u_{\mathrm{h}}) \, \phi_{j} \, \phi_{i} \dx{x}.
\end{gather*}
In total, we consider the following three systems for the standard Galerkin approach using Picard iterations:
\begin{subequations}
\begin{align}
    \nu \, \mat{K} \, \vec{u}^{n+1} + \mat{M}(\tilde{c},\vec{u}^{n}) \, \vec{u}^{n+1} &= \vec{d}, \label{eqn:superconductivity_a} \\
    \nu \, \mat{K} \, \vec{u}^{n+1} + \mat{M} \, \vec{u}^{n+1} + \vec{M}(\tilde{c},\vec{u}^{n}) &= \vec{d}, \label{eqn:superconductivity_b} \\
     \nu \, \mat{K} \, \vec{u}^{n+1} + \vec{M}(c,\vec{u}^{n}) &= \vec{d} \label{eqn:superconductivity_c}
\end{align}
\end{subequations}
with the respective nonlinearity $ \tilde{c} $.

Each way of decomposing the nonlinear function $ c $ results in different discrete systems with different unknowns. This influences the choice of the approximation space $ W_{\mathrm{h}} $, when applying the extended group finite element formulation. We consider piece-wise quadratic functions ($ \pp_{2} $) for \eqref{eqn:superconductivity_a} and piece-wise cubic functions ($ \pp_{3} $) for \eqref{eqn:superconductivity_b} and \eqref{eqn:superconductivity_c} as well as a quadrature rule that is exact for quartic polynomials ($ \ii_{4} $) for all systems. The group finite element method chooses $ W_{\mathrm{h}} = V_{\mathrm{h}} $ in each situation. The discrete problems solved by the (extended) group finite element method are given through
\begin{subequations}
\begin{align}
\nu \, \mat{K} \, \vec{u}^{n+1} + \tens{M}_{\mathrm{\tilde{c}}} : (\vec{u}^{n+1} \otimes \vec{\tilde{c}}^{n}) &= \vec{d}, \\
    \nu \, \mat{K} \, \vec{u}^{n+1} + \mat{M} \, \vec{u}^{n+1} + \mat{M}^{\mathrm{\tilde{c}}} \, \vec{\tilde{c}}^{n} &= \vec{d}, \\
    \nu \, \mat{K} \, \vec{u}^{n+1} + \mat{M}^{\mathrm{c}} \, \vec{c}^{n} &= \vec{d}
\end{align}
\end{subequations}
with appropriate point-wise definitions for the auxilliary variable $ \vec{c} $ or $ \vec{\tilde{c}} $.

\begin{figure} \centering
	\includegraphics{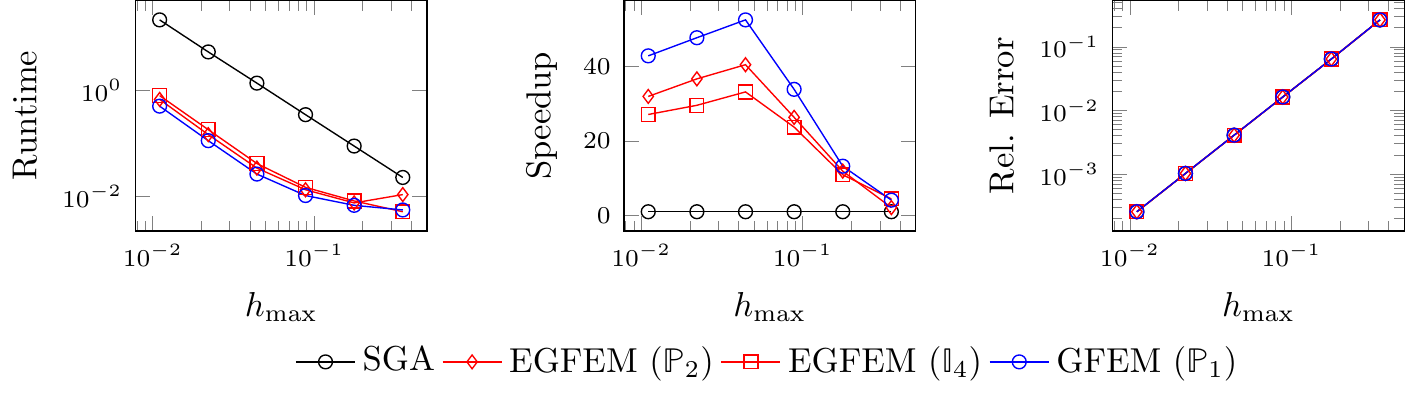}

    \caption{An overview of the performance of each approach using \eqref{eqn:superconductivity_a} with $ \nu = 1 $. The source term is assembled using an adaptive quadrature rule, which reduces the numerical errors from the numerical integration as much as possible. \textsl{From left to right:} Runtime in seconds to solve the problem, speedup associated to each formulation (with respect to the standard approach) and the relative error with respect to the exact solution using the approximate $L^{2}$-norm.}
    \label{fig:superconductivity}
\end{figure}

For the numerical tests, we consider $ \nu \in \{1,10^{-2},10^{-3}\} $, since for decreasing values of $ \nu $ the nonlinear term becomes more critical. A complete overview of these results for a fixed mesh can be seen in Table \ref{tbl:superconductivity}. The empty entries of the table denote the formulations, where the fixed point iteration either diverges or oscillates. For the situations in which the method converges, we see a noticeable decrease in the computational effort for the reformulated problems. Figure \ref{fig:superconductivity} depicts the performance of the discussed methods for $ \nu = 1 $ and the setting \eqref{eqn:superconductivity_a}. Here, we see the improvements in performance for the systems that avoid repeated assembly of the nonlinear forms. The error in Figure \ref{fig:superconductivity} is dominated by the approximation error, which leads to only a slight variation between the reformulated systems. The error introduced by approximating the nonlinearity with the GFEM in comparison to the exact reformulation of the EGFEM is negligible. In summary, a careful balance between the size of the reformulated system and the approximation errors must be found. In situations such as this example, Case~2 in Section \ref{subsec:egfem_choice} -- only approximating the nonlinearity -- is not only viable but also more efficient than the exact alternatives. In this case, the increased computational effort for exactly reformulating the nonlinear terms is unnecessary and can be avoided. This is most evident in the speedup, where the size of the reformulated problems directly dictates the methods efficiency. The group finite element method, which uses linear finite elements for $ W_{h} $, results in a much smaller problem, when compared to the quadratic finite elements or high-degree quadrature rule of the extended group finite element method, and as such has faster solution times.

\begin{table} \centering \renewcommand{\arraystretch}{1.15} \small
    \begin{tabular}{l | c | c | c | c | c | c | c | c}
        & \textbf{Model} & \textbf{Size} & \textbf{Time} & \textbf{Speedup} & \textbf{Time} & \textbf{Speedup} & \textbf{Time} & \textbf{Speedup} \\ \hline\hline
        \multirow{4}{*}{\STAB{\rotatebox[origin=c]{90}{\eqref{eqn:superconductivity_a}}}} & SGA & 4225  & 5.393 & 1.000 & 20.792 & \hphantom{0}1.000 & 42.256 & \hphantom{0}1.000 \\
        & EGFEM ($ \pp_2 $) & 20866 & 0.146 & 36.829 & \hphantom{0}0.577 & 36.025 & \hphantom{0}1.185 & 35.651 \\
        & EGFEM ($ \ii_4 $) & 53377 & 0.182 & 29.699 & \hphantom{0}0.702 & 29.627 & \hphantom{0}1.425 & 29.663 \\
        & GFEM              & 8450  & 0.112 & 47.941 & \hphantom{0}0.437 & 47.550 & \hphantom{0}0.894 & 47.289 \\ \hline \hline
        \multirow{4}{*}{\STAB{\rotatebox[origin=c]{90}{\eqref{eqn:superconductivity_b}}}} & SGA & 4225  & 4.651 & \hphantom{0}1.000 & 34.759 & \hphantom{0}1.000 & - & - \\
        & EGFEM ($ \pp_3 $) & 41474 & 0.076 & 61.189 & \hphantom{0}0.599 & 58.061 & - & - \\
        & EGFEM ($ \ii_4 $) & 53377 & 0.077 & 60.159 & \hphantom{0}0.579 & 60.009 & - & - \\
        & GFEM              & 8450  & 0.071 & 65.585 & \hphantom{0}0.534 & 65.037 & - & - \\ \hline \hline
        \multirow{4}{*}{\STAB{\rotatebox[origin=c]{90}{\eqref{eqn:superconductivity_c}}}} & SGA & 4225  & 7.295 & \hphantom{0}1.000 & - & - & - & - \\
        & EGFEM ($ \pp_3 $) & 41474 & 0.084 & 86.704 & - & - & - & - \\
        & EGFEM ($ \ii_4 $) & 53377 & 0.086 & 84.482 & - & - & - & - \\
        & GFEM              & 8450  & 0.079 & 92.769 & - & - & - & - \\ \hline
        \multicolumn{3}{c|}{} & \multicolumn{2}{c|}{$ \nu = 1 $} & \multicolumn{2}{c|}{$ \nu = 10^{-2} $} & \multicolumn{2}{c}{$ \nu = 10^{-3} $}
    \end{tabular}
    \caption{Comparison of the runtimes (in seconds) and speedup with respect to the standard Galerkin approach (SGA) for the different solution methods, when solving the simplified superconductivity model equations with varying values of $ \nu $. The empty fields (marked with -) designate formulations for which the fixed point iteration does not converge to the expected solution.}
    \label{tbl:superconductivity}
\end{table}

\subsection{Biochemical Reaction} \label{subsec:biochemical}

Let us consider the following problem
\begin{align*}
    -\Delta u(x) + \sigma \, \frac{u(x)}{k + u(x)} &= d(x), & x &\in \Omega \mdef [0,1] \times [0,1], \\ 
    u(x) &= x_{1} \, x_{2} \, (x_{1} + x_{2}), & x &\in \Gamma,
\end{align*}
where the reaction parameters $ \sigma $ and $ k $ as well as the source term $ d $ are positive. For simplicity, we choose $ \sigma = 1 = k $ and determine $ d $ with the method of manufactured solutions. Similar to Section \ref{subsec:superconductivity}, we have a situation, where the nonlinear term $ c $ can either be written as
\begin{align*}
    \int_{\Omega_{\mathrm{h}}} c(u_{\mathrm{h}},\nabla u_{\mathrm{h}}) \, \phi_{i} \dx{x} &= \mat{M}(\tilde{c},\vec{u}) \, \vec{u}, & \tilde{c}(u,\nabla u) &= \frac{\sigma}{k + u}, \\
     \int_{\Omega_{\mathrm{h}}} c(u_{h},\nabla u_{\mathrm{h}}) \, \phi_{i} \dx{x} &= \vec{M}(c,\vec{u}). 
\end{align*}
However, in light of the previous results, we only consider the first case, which delivers the following system for the SGA with Picard iterations:
\begin{gather*}
	\mat{K} \, \vec{u}^{n+1} + \mat{M}(\tilde{c},\vec{u}^{n}) \, \vec{u}^{n+1} = \vec{d}.
\end{gather*}%
Here, we have a situation, where the nonlinearity -- a rational function in terms of $ u $ -- does not have a matching finite element. Therefore, we consider Case~2 and Case~3 of Section \ref{subsec:egfem_choice}. Using the (extended) group finite element formulation, the system can be rewritten as
\begin{gather*}
    \mat{K} \, \vec{u}^{n+1} + \tens{M}_{\mathrm{\tilde{c}}} : (\vec{u}^{n+1} \otimes \vec{\tilde{c}}^{n}) = \vec{d}, \quad 
        (\vec{\tilde{c}}^{n+1})_{i} = \frac{\sigma}{k + u_{\mathrm{h}}^{n+1}(x_{i}^{\mathrm{c}})}
\end{gather*}
with $ x_{i}^{\mathrm{c}} $ denoting the degrees of freedom for $ \tilde{c}_{\mathrm{h}} $. 

\begin{figure} \centering
	\includegraphics{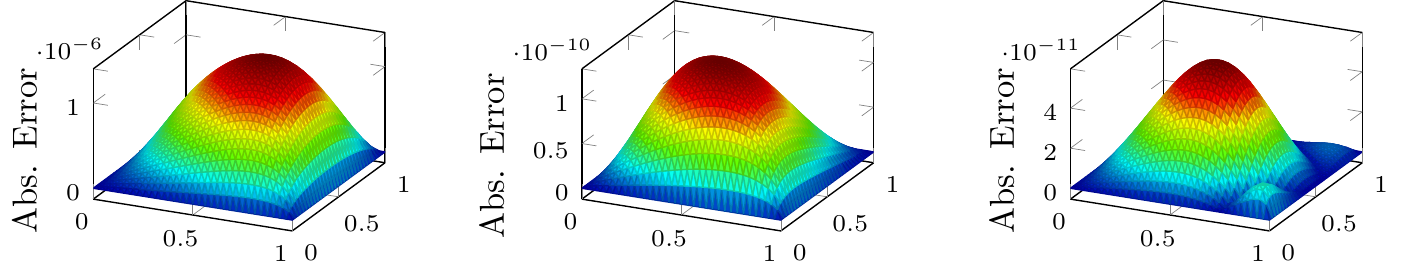}

    \caption{Convergence of quadrature rule. \textsl{From left to right:} absolute error between the solutions using first and second order, second and third order and third and fourth order quadrature rules, respectively. \label{fig:biochemical_convergence}}
\end{figure}

The only problem that remains is which quadrature rule makes the most sense in order to evaluate the nonlinear term. If we look at the original term and use the fact that $ u_{\mathrm{h}}|_{K} \in \pp_{1} $, then we can see that
\begin{align*}
    \int_{\Omega_{\mathrm{h}}} \frac{\sigma \, \phi_{j}(x) \, \phi_{i}(x)}{k + u_{\mathrm{h}}(x)} \dx{x} &= \int_{\Omega_{\mathrm{h}}} \frac{P(x)}{Q(x)} \dx{x},
\end{align*}
where $ P $, $ Q $ are at most polynomials of degree 2 and 1, respectively. Using this knowledge, it makes sense that the rational term $ P / Q $ behaves at most quadratically, and, therefore, a quadrature rule, which is exact for quadratic polynomials, delivers sufficient results. This hypothesis is confirmed in Figure \ref{fig:biochemical_convergence}, where quadrature rules of increasing order are compared.

\begin{figure} \centering
	\includegraphics{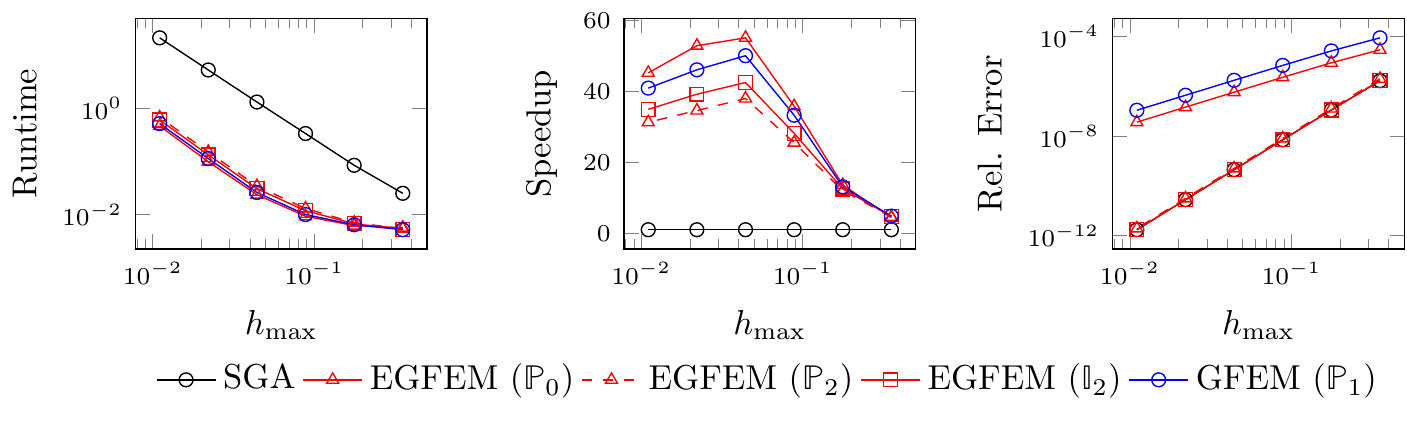}

    \caption{Comparison of numerical results for the biochemical equation. \textsl{From left to right:} average time to compute the solution, speedup associated to each approach and the relative error with respect to the manufactured solution using the approximate $ L^{2} $-norm. \label{fig:biochemical}}
\end{figure}

The results shown in Figure \ref{fig:biochemical} lead to interesting conclusions. First, we see that the quadrature elements ensure that the standard approach is conserved. In contrast, the approximations, while not being exact reformulations of the original problem, offer excellent results even with piece-wise constant elements. The errors in Figure \ref{fig:biochemical} show the influence of the reformulation technique. Interestingly, the higher order approximation using quadratic finite elements displays the same accuracy as the SGA and equivalent reformulation using the quadrature rule. Inexplicably, the lower order approximation using piecewise constant functions ($ \pp_{0} $) performs better than the original group finite element method using piecewise linear functions ($ \pp_{1} $). Perhaps this result is similar to reports that the GFEM converges with a higher order than the SGA in certain norms, see \cite{Christie1981,Fletcher1983}. In addition to the superior accuracy, the EGFEM using $ \pp_{0} $ elements computes the solution faster than the GFEM, which can be explained by the smaller problem size.

\subsection{$ p $-Laplace Equation} \label{subsec:plaplace}

The $ p $-Laplacian, which generalizes the Laplacian operator, defines a quasilinear elliptic partial differential operator. The $ p $-Laplace equation with Dirichlet boundary conditions is given through
\begin{align*}
	-\nabla \cdot \left(\| \nabla u(x) \|^{p-2} \, \nabla u(x)\right) &= d(x), & x &\in \Omega, \\
	u(x) &= u_{\mathrm{D}}(x), & x &\in \Gamma.
\end{align*}
For the case $ d \equiv 1 $, $ u_{\mathrm{D}} \equiv 0 $ and $ \Omega = \{x \in \rr^{2} : \| x \| \leq 1\} $, an explicit solution is given by
\begin{gather*}
	u(x) = 2^{-\frac{1}{p-1}} \, \frac{p-1}{p} \left(1 - \| x \|^{\frac{p}{p-1}}\right)
\end{gather*}
for $ p \in (1,\infty) $, see, for example, \cite{Kawohl1990}. In our setting of \eqref{eqn:nonlinear_pde}, we have
\begin{gather*}
	a(x,u,\nabla u) = \|\nabla u\|^{p-2}, \quad b(u) = u, \quad c(x,u,\nabla u) = 0, \quad d(x) = 1
\end{gather*}
with $ \Gamma_{\mathrm{D}} = \Gamma $ and $ \Gamma_{\mathrm{N}} = \emptyset $. The discretization using linear finite elements leads to the Picard iterations for the standard approach
\begin{gather*}
	\mat{K}(a,\vec{u}^{n}) \, \vec{u}^{n+1} = \vec{d}
\end{gather*}
and for the extended group finite element formulation
\begin{gather*}
	\tens{K}_{\mathrm{a}} : (\vec{u}^{n+1} \otimes \vec{a}^{n}) = \vec{d}, \quad \vec{a}^{n+1} = a(\Pi_{\nabla \mathrm{u}}^{\mathrm{a}} \cdot_{2} \vec{u}^{n+1}).
\end{gather*}
Because the gradient $ \nabla u_{\mathrm{h}} $ is discontinuous along the interfaces of the elements, it is not possible to interpolate the nonlinearity $ a $ onto the trial space $ V_{\mathrm{h}} $. Hence, the GFEM is not applicable. For this example, we consider the EGFEM with $ W_{\mathrm{h}} $ using piece-wise constant elements ($ \pp_{0} $) and the midpoint quadrature rule ($ \ii_{1} $), which is exact for linear functions.

\begin{figure} \centering
	\includegraphics{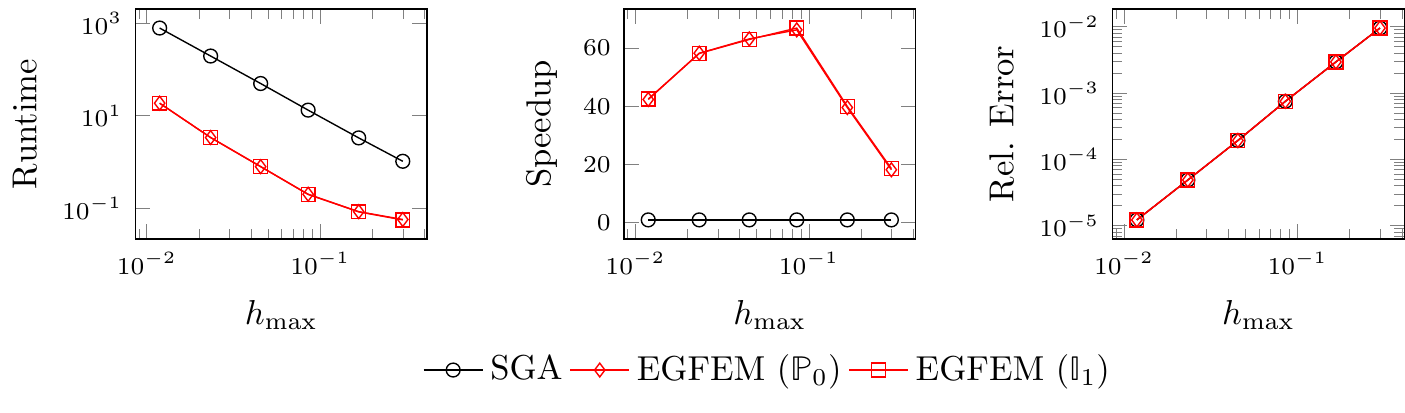}

    \caption{An overview of the performance of each approach for the $ p $-Laplace equation using a regular mesh that is successively refined through splitting. \textsl{From left to right:} The average amount of time to solve the problem for increasingly finer meshes, the speedup describes the ratio between the reference runtime (SGA) and the considered approach and the relative error with respect to the $ L^{2} $-norm from the exact solution.}
    \label{fig:pLaplace}
\end{figure}

The computational advantages of the extended group finite element formulation are depicted in Figure \ref{fig:pLaplace} for $ p = 3/2 $. Additionally, the extended formulations are equivalent with the standard Galerkin approach. In such situations, performance gains are achieved at the cost of increasing storage requirements, e.g., for storing the third-order tensors, which can be very large but are also sparse, and the intermediate variables, which store the information about the nonlinearity.

\subsection{Minimal Surface Equation} \label{subsec:minimal}

For this example, we consider the minimal surface equation, which is a quasi-linear partial differential equation, i.e.,
\begin{align*}
    -\nabla \cdot \left(\frac{\nabla u(x)}{\sqrt{1 + \| \nabla u(x) \|^{2}}}\right) &= d(x), & x &\in \Omega \,\mdef [0,1] \times [0,1], \\
    u(x) &= u_{\mathrm{D}}(x) = x_{1} \, x_{2} \, (x_{1} + x_{2}), & x &\in \Gamma.
\end{align*}
If the source $ d $ vanishes, then the solution describes the function whose graph has the minimal area of all graphs over $ \Omega $ with the fixed values on the boundary \cite{Farina2007}. For the numerical investigations, the source term $ d $ is chosen according to the method of manufactured solutions \cite{Roache2019}, where the desired solution is given through the extension of $ u_{\mathrm{D}} $ onto the entire domain $ \bar{\Omega} $. This allows us to study the accuracy of the method by comparing the numerical solutions to the exact solution, while also comparing the performance of the extended group finite element formulation with the standard Galerkin approach. In relation to \eqref{eqn:nonlinear_pde}, this example is given through
\begin{gather*}
    a(x,u,\nabla u) = (1 + \| \nabla u \|^{2})^{-\frac{1}{2}}, \quad b(u) = u, \quad 
    c \equiv 0, \quad d(x) = -\nabla \cdot \left(\frac{\nabla u_{\mathrm{D}}(x)}{\sqrt{1 + \| \nabla u_{\mathrm{D}}(x) \|^{2}}}\right)
\end{gather*}
with $ \Gamma_{\mathrm{D}} = \Gamma $ and $ \Gamma_{\mathrm{N}} = \emptyset $. After discretizing with linear Lagrange elements, the standard Galerkin approach delivers the following nonlinear discrete problem:
\begin{gather*}
    \mat{K}(a,\vec{u}^{n}) \, \vec{u}^{n+1} = \vec{d}.
\end{gather*}
In contrast, the extended group finite element method solves
\begin{gather*}
    \tens{K}_{\mathrm{a}} : (\vec{u}^{n+1} \otimes \vec{a}^{n}) = \vec{d}, \quad \vec{a}^{n+1} = a(\Pi_{\nabla \mathrm{u}}^{\mathrm{a}} \cdot_{2} \vec{u}^{n+1}).
\end{gather*}
Once again, the GFEM is not viable, since the gradient of the linear approximation $ u_{\mathrm{h}} $ is not defined on the nodes. Therefore, the extension presented in this work is necessary. In light of the fact that $ a $ is constant on each element, we consider the cases that $ W_{\mathrm{h}} $ is given through piece-wise constant functions ($ \pp_{0} $) and using a quadrature rule that is exact for linear and, therefore, constant functions ($ \ii_{1} $).

\begin{figure} \centering
	\includegraphics{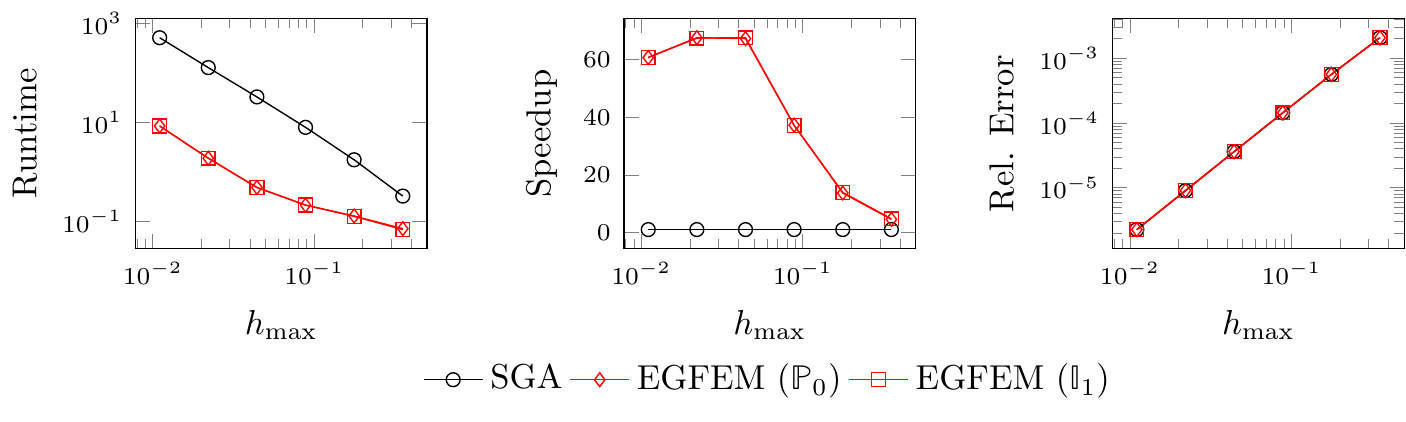}

    \caption{An overview of the performance of each approach for the minimal surface equation using a regular mesh that is successively refined through splitting. \textsl{From left to right:} The average amount of time to solve the problem for increasingly finer meshes, the speedup describes the ratio between the reference runtime (SGA) and the considered approach, and the relative error with respect to the approximate $ L^{2} $-norm from the exact solution.}
    \label{fig:minsurface}
\end{figure}

Figure \ref{fig:minsurface} illustrates the advantages of the EGFEM. The average runtime and, consequently, the respective speedup testify to the performance increases by eliminating the repeated assembly of the nonlinear problem. For this benchmark problem, the reformulation using intermediate variables for the nonlinearity is numerically exact as can be seen in convergence of the absolute error of the solutions under mesh refinement. In summary, the extension to more general approximation spaces $ W_{\mathrm{h}} $ is necessary in order to apply the group finite element formulation to this problem, which results in a numerically equivalent \emph{structured} problem that can be solved faster than the standard approach.

\subsection{Discussion}

In the course of the numerical investigations, a collection of observations have arisen that may serve as a set of guidelines, when choosing the approximation space for the extended group finite element method. The group finite element method offers a natural approximation of nonlinear terms in order to improve the numerical performance, as shown in \cite{Christie1981,Fletcher1983}. However, there are situations where the GFEM fails. For example, the $ p $-Laplace equation in Section \ref{subsec:plaplace} and the minimal surface equation in Section \ref{subsec:minimal} break down, when linear Lagrange elements are used for $ V_{\mathrm{h}} $. In this case, the generalizations of our extended group finite element method play a key role.

By introducing a general approximation space $ W_{\mathrm{h}} $ for the nonlinearity, the extended group finite element method becomes a flexible and powerful tool, that offers a problem-tailored approach. However, it is difficult to formulate concrete guidelines, since the problem and application \emph{both} play an important role in the choice of the approximation. In Section \ref{subsec:egfem_choice}, we defined three cases for the choice of $ W_{\mathrm{h}} $ with respect to the nonlinearity. While Case~1 exactly reformulates the nonlinearity and Case~3 the original problem, Case~2 explicitly chooses an approximation of the considered function. Case~1 requires that the chosen trial space $ V_{\mathrm{h}} $ and nonlinear function result in a class of functions that can be exactly described by an existing finite element, which may not always be possible. In comparison, Case~3, which embeds a quadrature rule into an approximation space, conserves the original discrete problem but often leads to a large number of additional variables. In between these two cases, Case~2 is perhaps the most important. First, the GFEM can be recovered by choosing $ W_{\mathrm{h}} = V_{\mathrm{h}} $. However, it is also possible to choose a lower or higher order approximation based on the application and problem. Therefore, Case~2 of the EGFEM allows for much more control in comparison to the GFEM.

The numerical examples presented throughout this section clearly show that the construction of $ W_{\mathrm{h}} $ crucially depends on the problem and application. For problems where the approximation error through $ V_{\mathrm{h}} $ dominates, Case~2 is often the most viable. For example, the low-order approximation with the GFEM outperforms the considered EGFEM systems in Section \ref{subsec:superconductivity}. In this situation, an exact reformulation or high-order approximation is unnecessary and can be avoided. However, if the approximation error through $ V_{\mathrm{h}} $ is negligible, then the contrary is true. Then, the choice for the approximation space $ W_{\mathrm{h}} $ depends on the application. For example, the setting in Section \ref{subsec:biochemical} leads to a situation, where the performance of the different approaches varies considerably. While the approximations using quadratic elements or a high-order quadrature rule outperform the SGA, they also take longer to solve than the simpler approximations using constant or linear elements. However, the relative error remains larger for the simpler approximations than the higher-order approximations. Therefore, a balance has to be found between the quality of the approximation and the time required to solve the resulting problem.

\section{Conclusion}

The assembly of the nonlinear systems accounts for a large amount of computational overhead. Therefore, techniques like the group finite element formulation are crucial for time-sensitive computations. An important and relevant example is the optimal control of a nonlinear PDE system. Because gradient-based methods require the repeated solution of the forward and backward problems, the original and our extended group finite element methods become invaluable, since large portions of the systems only need to be computed once. However, the original group finite element formulation has its limitations. By interpolating all nonlinear terms onto the trial space $ V_{\mathrm{h}} $, the resulting system is only an approximation of the original problem. At the same time, many difficulties arise for nonlinear problems that depend on gradient information. In such situations, the original group finite element formulation may not even be applicable. In order to allow for more flexibility as well as improve the accuracy of the reformulated systems, we extend the group finite element method to use general finite element approximation spaces. As seen in Section \ref{sec:results}, this allows for the exact reformulation of the original problem in many cases as well as the handling of nonlinear functions that depend on the gradient of the approximate solution $ u_{\mathrm{h}} $. The additional extension of the group finite element formulation to multi-linear forms using tensors introduces structure to problems that can be used even without the introduction of intermediate variables, see, for example, Sections \ref{subsec:quadratic} and \ref{subsec:burgers}. These adaptations lead to a more efficient approximation of nonlinear problems with similar computational benefits as the group finite element method. While our results use the Picard iterations, it is reasonable to expect qualitatively similar results when using Newton's method. This is because the most important aspect of the extended group finite element method is the elimination of the computational overhead associated with repeated assembly. A real-world application using our proposed method is presented in \cite{TolleICIAM}. Additionally, further extensions including the combination with model order reduction techniques, similar to \cite{Dickinson2010,Wang2015}, are currently being investigated. Here, the discrete empirical interpolation method (DEIM) plays an important role in the efficient evaluation of the nonlinear terms, see \cite{Chaturantabut2010,Antil2014}. In standard finite element applications, alternative formulations of the discrete empirical interpolation method are necessary in order to accommodate for the influence of neighboring nodes, as described in \cite{Tiso2013}. However, an interesting aspect of our approach that plays to the original advantages of DEIM is the point-wise evaluation of the nonlinearity.

\section*{Acknowledgments}

We would like to acknowledge the support of the BMBF \textsc{proMT}-project as well as the BMWi MathEnergy project.

\bibliographystyle{siam}
\bibliography{references.bib}

\end{document}